\documentclass[12pt,a4paper]{article}
\usepackage{amssymb}
\usepackage{epsfig}
\usepackage{amsmath}
\usepackage{amsfonts}
\usepackage{pgf,tikz}
\usepackage{enumerate}
\usepackage{cite}
\usepackage{stmaryrd}
\usepackage[autobold]{mathfixs}
\usepackage[colorlinks=true, linkcolor=black, citecolor=black, urlcolor=black]{hyperref}

\usepackage[vmargin=2cm, hmargin=2cm]{geometry}

\newcommand{\R}{{{\Bbb R}}}
\DeclareMathOperator{\Id}{Id}
\DeclareMathOperator{\dif}{d}
\newcommand{\KS}{{}^{(KS)}\!\!\!}

\newtheorem{theorem}{\sc Theorem}[section]
\newtheorem{proposition}[theorem]{\sc Proposition}
\newtheorem{lemma}[theorem]{\sc Lemma}
\newtheorem{definition}[theorem]{\sc Definition}
\newtheorem{remark}[theorem]{\sc Remark}
\newtheorem{corollary}[theorem]{\sc Corollary}
\newtheorem{example}[theorem]{\sc Example}
\def\qed{\hbox to 0pt{}\hfill$\rlap{$\sqcap$}\sqcup$\medbreak}

\title{Existence and uniqueness of solution for Stieltjes differential equations with several derivators}

\author{Ignacio M\'arquez Alb\'es and F. Adri\'an F. Tojo}

\date{}
\begin{document}
 \maketitle

\begin{center}  {\small Departamento de Estat\'{\i}stica, An\'alise Matem\'atica e Optimizaci\'on \\ Instituto de Matem\'aticas \\ Universidade de Santiago de Compostela \\ 15782, Facultade de Matem\'aticas, Campus Vida, Santiago, Spain. \\ e-mail: ignacio.marquez@usc.es, fernandoadrian.fernandez@usc.es}
\end{center}

\medbreak

\begin{abstract}
In this paper, we study some existence and uniqueness results for systems of differential equations in which each of equations of the system involves a different Stieltjes derivative. Specifically, we show that this problems can only have one solution under the Osgood condition, or even, the Montel--Tonelli condition. We also explore some results guaranteeing the existence of solution under these conditions. Along the way, we obtain some interesting properties for the Lebesgue--Stieltjes integral associated to a finite sum of nondecreasing and left--continuous maps, as well as a characterization of the pseudometric topologies defined by this type of maps.
\end{abstract}

\medbreak


\medbreak


\section{Introduction}

Stieltjes differential equations have gained popularity in the recent years. The main difference with respect to regular differential problems is the presence of the Stieltjes derivative, a modification of the usual derivative on the real line through a nondecreasing and left--continuous map. This change allows us to study impulsive differential equations and equations on time scales in a unified framework, see for example \cite{LoRo14,FriLo17}, or even \cite{FeMesSla12} for the corresponding integral formulation counterpart.

The usual setting for Stieltjes differential equations in the literature involves a single derivator either in its theoretical --see for example \cite{FriLo17,LoMa18,LoMa20,MaMon20,MonSat17,FerTo20,SatSmy20}-- or numerical studies \cite{Fernandez2020,FerTo20}. This is also the case for other differential problems involving Stieltjes derivatives such as in \cite{LoMaRo20,MonSat19,SatSmyr20}, or even the corresponding integral counterparts. Nevertheless, it is the new setting of differential problems with Stieltjes derivatives that offers the possibility of a new type of problems: systems of differential equations in which each of the components is differentiated with respect to a different nondecreasing and left--continuous function. This was the case in papers such as \cite{LoMa19,LoMaMon18}, or even \cite{FriTo20}, where the authors considered differentiation with respect to functions that are not necessarily monotonous. Here, we aim to improve the work along this line regarding systems of equations. Specifically, we will consider  maps $g_i:\mathbb R\to\mathbb R^n$, $i=1,2,\dots,n$, such that each $g_i$ is nondecreasing and left--continuous, and we will discuss some existence and uniqueness results for the system
\begin{equation}\label{ivpintro}
x'_{g_i}(t)=f_i(t,x(t)),\quad x_i(t_0)=x_{0,i},\quad i=1,2,\dots,n,
\end{equation}
where $x'_{g_i}$ denotes the Stieltjes derivative of $x$ with respect to $g_i$ in the sense presented in \cite{LoRo14}. In this setting, we build on the work in \cite{LoMa19}, adapting some of the results in \cite{FriLo17,MaMon20} to the context of Stieltjes differential equations with several derivators.

The paper is structured as follows: first, in Section~\ref{sectionintegral}, we present the basic tools for the Lebesgue--Stieltjes integration on the real line defined in terms of a nondecreasing and left--continuous map. We also obtain some interesting results regarding the Lebesgue--Stieltjes outer measure, as well as a fundamental property regarding the Lebesgue--Stieltjes measure for the finite sum of nondecreasing and left--continuous functions --more information on the Lebesgue-Stieltjes measure can be found in \cite{Winter1997,FriLo17,LoRo14,AthLa06,Burk07,Hil63} and for the more general Kurzweil-Stieltjes integral in \cite{MonSlaTvr18}. Next, in Section~\ref{sectionderivative}, we introduce the Stieltjes derivative in the sense of \cite{FriLo17} and we explore some concepts of continuity in a similar fashion to \cite{FriLo17,LoMa19}. In particular, we discuss some of the limitations of the mentioned work. Furthermore, throughout this section we obtain important information regarding the pseudometric topology defined by a nondecreasing and left--continuous function, showing that it can be fully characterized in terms of some interesting sets related to such map. Finally, in Section~\ref{sectionivp}, we turn our attention to the study of problems of the form \eqref{ivpintro}. First, we continue the study of everywhere solutions started in \cite{LoMa19} and, later, we follow the arguments in \cite{MaMon20} to obtain some existence and uniqueness results involving Osgood and Montel--Tonelli conditions.

\section{The Lebesgue--Stieltjes measure}\label{sectionintegral}
Throughout this paper, we will make us of the Lebesgue--Stieltjes integral associated to nondecreasing and left--continuous functions. This integral is constructed as the integral with respect to a measure defined in terms of the mentioned map through the classical Carath\'eodory's extension theorem, see for example \cite{AthLa06,Burk07,Mun53,Sche97,Ru87}. Specifically, given a nondecreasing and left--continuous map $g:\mathbb R\to\mathbb R$, and denoting by $\mathcal P(\mathbb R)$ the set of all subsets of $\mathbb R$, we define the map $\mu_g^*:\mathcal P(\mathbb R)\to[0,+\infty]$ as
\begin{equation}\label{mugextintab}
\mu_g^*(A)=\inf\left\{\sum_{n=1}^\infty (g(b_n)-g(a_n)): A\subset \bigcup_{n=1}^\infty [a_n,b_n),\ \{[a_n,b_n)\}_{n=1}^\infty\subset \mathcal C\right\},
\end{equation}
with  $\mathcal C=\{[a,b):a,b\in\mathbb R,\ a<b\}$. The map $\mu_g^*$ is an outer measure and, by considering its restriction to the following $\sigma$--algebra,
\[\mathcal{LS}_g=\{A\in\mathcal P(\mathbb R)\, : \, \mu_g^*(E)=\mu_g^*(E\cap A)+\mu_g^*(E\setminus A)\mbox{ for all }E\in \mathcal P(\mathbb R)\},\]
we obtain the Lebesgue--Stieltjes measure associated to $g$, which we denote by $\mu_g$.
\begin{remark}\label{measureinterv}
	Every Borel set belongs to $\mathcal{LS}_g$. In particular, this means that $\mathcal C\subset \mathcal{LS}_g$. Furthermore, we have that $\mu_g([a,b))=g(b)-g(a)$ for any $[a,b)\in\mathcal C$.
\end{remark}

For simplicity, in what follows we will use the term ``$g$--measurable'' for a set or function to refer to $\mu_g$--measurability in the corresponding sense; and we will denote the integration with respect to $\mu_g$ as
\[\int_X f(s)\dif g(s).\]
In a similar way, we will replace $\mu_g$ by $g$ in other expressions such as ``$P$ holds for $\mu_g$--a.a. $x\in X$'' or  ``$P$ holds $\mu_g$--a.e. in $X$''. Along these lines, it is important to note that the set
\[C_g:=\{ t \in \mathbb R \, : \, \mbox{$g$ is constant on $(t-\varepsilon,t+\varepsilon)$ for some $\varepsilon>0$} \},\]
i.e. the set of points around which $g$ is constant, has null $g$--measure, as pointed out in \cite[Proposition~2.5]{LoRo14}. Furthermore, observe that, by definition, $C_g$ is open.

The aim of this section is to show that the expression used to compute $\mu_g$, \eqref{mugextintab}, can be simplified, as well as to prove some interesting properties regarding the Lebesgue--Stieltjes measure associated to the sum of a finite family of nondecreasing and left--continuous functions.

We begin by showing that \eqref{mugextintab} can be simplified. Specifically, we will show that the infimum in that expression can be considered over an smaller set, namely, assuming that the families in $\mathcal C$ are pairwise disjoint. To that end, we introduce the following lemma.
\begin{lemma}\label{lemrecubr}
	Let $g:\mathbb R\to\mathbb R$ be a nondecreasing and left--continuous function. For every $\{[a_n,b_n)\}_{n\in\mathbb N}\subset\mathcal C$, there exists $\mathcal V=\{[c_n,d_n)\}_{n\in\mathbb N}\subset\mathcal C$ such that the sets in $\mathcal V$ are pairwise disjoint and
	\[\bigcup_{n\in\mathbb N}[a_n,b_n)=\bigcup_{n\in\mathbb N}[c_n,d_n),\quad\quad\quad \sum_{n\in\mathbb N}(g(d_n)-g(c_n))\le \sum_{n\in\mathbb N} (g(b_n)-g(a_n)).
	\]
\end{lemma}
\noindent
{\bf Proof.}
	Let $U=\bigcup_{n\in\mathbb N}[a_n,b_n)$ and $C_U$ be the set of all connected components of $U$. First, note that the set $C_U$ is at most countable.
	Secondly, observe that all the elements of $C_U$ are connected subsets of $\mathbb R$. Thus, we have that they are intervals (including the whole $\mathbb R$) or singletons. Nevertheless, observe that an element of $C_U$ cannot be a singleton as each point of $U$ belongs to $[a_n,b_n)$ for some $n_0\in\mathbb N$. 
	Furthermore, we claim that $\sup I\not\in I$ for any $I\in C_U$ bounded from above. Indeed, let $I\in C_U$ be bounded from above and suppose that $\sup I\in I\subset U$. In this conditions, there exists $n_1\in\mathbb N$ such that $\sup I\in[a_{n_1},b_{n_1})$. Hence, we have that the set $I\cup[a_{n_1},b_{n_1})$ is a connected set containing $I$, which is a contradiction with $I\in C_U$.
	Therefore, $C_U$ is, by construction, an at most countable collection of pairwise disjoint sets of the form $(a,b)$, $[a,b)$, $[a,+\infty)$,  $(-\infty, b)$, $a,b\in\mathbb R$; or $C_U=\{\mathbb R\}$.

	For each $I\in C_U$ and define $\mathcal F_I\subset\mathcal C$ as follows:
	\begin{enumerate}
		\item[$\cdot$] if $I=(a,b)$, $a,b\in\mathbb R$, then  $\mathcal F_I=\displaystyle \left\{\left[a+(b-a)/(n+1),a+(b-a)/n\right)\right\}_{n\in\mathbb N}$;
		\item[$\cdot$] if $I=[a,b)$, $a,b\in\mathbb R$, then  $\mathcal F_I=\displaystyle \left\{\left[a+(n-1)(b-a)/n,a+n(b-a)/(n+1)\right)\right\}_{n\in\mathbb N}$;
		\item[$\cdot$] if $I=[a,+\infty)$, $a\in\mathbb R$, then $\mathcal F_I=\left\{\left[a+n-1,a+n\right)\right\}_{n\in\mathbb N}$;
		\item[$\cdot$] if $I=(-\infty,b)$, $b\in\mathbb R$, then $\mathcal F_I=\left\{\left[b-n,b-n+1\right)\right\}_{n\in\mathbb N}$,
		\item[$\cdot$] if $I=\mathbb R$, then $\mathcal F_I=\{[n,n+1)\}_{n\in\mathbb Z}$.
	\end{enumerate}
	Proceeding this way, for each $I\in C_U$ we find a countable pairwise disjoint family contained in $\mathcal C$, $\mathcal F_I$, such that $I=\bigcup_{J\in\mathcal F_I} J$. Furthermore, for each $I\in C_U$, it follows from Remark \ref{measureinterv} and the fact that $\mu_g$ is a measure that
	\begin{equation}\label{telescopingphi}
	\sum_{J\in\mathcal F_I} (g(\sup J)-g(\inf J))=\sum_{J\in\mathcal F_I} \mu_g(J)=\mu_g\left(\bigcup_{J\in\mathcal F_I} J\right)=\mu_g(I).
	\end{equation}

	Define $\mathcal V=\bigcup_{I\in C_U}\mathcal F_I$. First, observe that $\mathcal V$ is a countable set by definition. Furthermore, $\mathcal V\subset \mathcal C$. Hence, we can write $\mathcal V=\{[c_n,d_n)\}_{n\in\mathbb N}$ for some $c_n,d_n\in~\mathbb R$. Let us show that $\mathcal V$ satisfies the properties in the statement of the result.

	First, note that the sets in $\mathcal V$ are pairwise disjoint. Indeed, let $[c_n,d_n)$ and  $[c_m,d_m)$ be two elements of $\mathcal V$. If they belong to the same connected component, $I\in C_U$, then, by construction of $\mathcal F_I$, we have that $[c_n,d_n)\cap [c_m,d_m)=\emptyset$. Otherwise, $[c_n,d_n)\in I$ and $[c_m,d_m)\in I'$ for some $I, I'\in C_U$, $I\not= I'$. Then, the definition of connected component guarantees that $I\cap I'=\emptyset$, which yields $[c_n,d_n)\cap [c_m,d_m)=\emptyset$. Hence, the family $\mathcal V$ is pairwise disjoint. Furthermore,
	\[U=\bigcup_{I\in C_U} I=\bigcup_{I\in C_U}\left(\bigcup_{J\in\mathcal F_I} J\right)=\bigcup_{V\in\mathcal V} V=\bigcup_{n\in\mathbb N}[c_n,d_n).\]
	Finally, using \eqref{telescopingphi} and Remark \ref{measureinterv}, we have that
	\begin{align*}
	\sum_{n\in\mathbb N}(g(d_n)-g(c_n))= & \sum_{I\in C_U}\left(\sum_{J\in\mathcal F_I} (g(\sup J)-g(\inf J))\right)\\ =&\sum_{I\in C_U}\mu_g(I)=\mu_g(U)  \le\sum_{n\in\mathbb N}(g(b_n)-g(a_n)),
	\end{align*}
	where the last inequality follows from \eqref{mugextintab}.
\qed

The following result contains a characterization of the outer measure $\mu_g^*$. The result follows from Lemma~\ref{lemrecubr} by considering the relations between the infima involved.
\begin{theorem}\label{teoroutermeasureg}
	Let $g:\mathbb R\to\mathbb R$ be a nondecreasing and left--continuous functions and $\mu_g^*$ be as in \eqref{mugextintab}. Then, for any $A\in\mathcal P(\mathbb R)$,
	\begin{equation}\label{mugextintab2}
	\mu^*_g(A)=\inf\left\{\sum_{n=1}^\infty (g(b_n)-g(a_n)): A\subset \bigcup_{n=1}^\infty [a_n,b_n),\ \{[a_n,b_n)\}_{n=1}^\infty\subset \mathcal C\mbox{ pairwise disjoint}\right\}.
	\end{equation}
\end{theorem}

Theorem \ref{teoroutermeasureg} not only provides an easier way to compute the outer measure of sets, but it is also a fundamental tool for the proof of the following result which, to the best of our knowledge, is not available in the existing literature on the topic of Lebesgue--Stieltjes measures. Essentially, Proposition~\ref{gsum} guarantees that, given a finite family of nondecreasing and left--continuous functions, the measure defined as the sum of the corresponding Lebesgue--Stieltjes measures is, in fact, the Lebesgue--Stieltjes measure associated to the corresponding sum of nondecreasing and left--continuous functions.

\begin{proposition}\label{gsum}
	Let $g_1,g_2,\dots,g_n:\mathbb R\to\mathbb R$ be a family of nondecreasing and left--continuous functions and define $\widehat g:\mathbb R\to\mathbb R$ as
	\begin{equation}\label{ghatlip}
		\widehat g(t)=g_1(t)+g_2(t)+\dots+g_n(t),\quad t\in\mathbb R.
	\end{equation}
	Then, for any $E\in\mathcal P(\mathbb R)$,
	\begin{equation}\label{gsumext}
	\mu_{\widehat g}^*(E)=\mu_{g_1}^*(E)+\mu_{g_2}^*(E)+\dots+\mu_{g_n}^*(E).
	\end{equation}
\end{proposition}
\noindent
{\bf Proof.} We shall only prove the result for $n=2$, as the general case can be deduced from this.

	Let $E\in\mathcal P(\mathbb R)$. Then, computing the corresponding outer measures as in \eqref{mugextintab}, we have that
	\begin{align*}
	\mu_{\widehat g}^*(E)&=\inf\left\{\sum_{n=1}^\infty (\widehat g(b_n)-\widehat g(a_n)): E\subset \bigcup_{n=1}^\infty [a_n,b_n)\right\}\\
	&=\inf\left\{\sum_{n=1}^\infty [(g_1(b_n)-g_1(a_n))+(g_2(b_n)-g_2(a_n))]: E\subset \bigcup_{n=1}^\infty [a_n,b_n)\right\}\\
	&=\inf\left\{\sum_{n=1}^\infty (g_1(b_n)-g_1(a_n))+\sum_{n=1}^\infty(g_2(b_n)-g_2(a_n)): E\subset \bigcup_{n=1}^\infty [a_n,b_n)\right\}\\
	&\ge \sum_{i=1}^2 \inf\left\{\sum_{n=1}^\infty (g_i(b_n)-g_i(a_n)): E\subset \bigcup_{n=1}^\infty [a_n,b_n)\right\}=\mu_{g_1}^*(E)+\mu_{g_2}^*(E).
	\end{align*}

	For the reverse inequality, let $\varepsilon>0$. It follows from \eqref{mugextintab2} that there exist $\mathcal R_1=\{[a_{1,n},b_{1,n})\}_{n=1}^\infty$ and $\mathcal R_2=\{[a_{2,m},b_{2,m})\}_{m=1}^\infty$, each of them pairwise disjoint, such that
	\begin{align*}
	&E\subset \bigcup_{n=1}^\infty [a_{1,n},b_{1,n}),&&\sum_{n=1}^\infty (g_1(b_{1,n})-g_1(a_{1,n}))\le \mu_{g_1}^*(E)+\frac{\varepsilon}{2},\\
	&E\subset \bigcup_{m=1}^\infty [a_{2,m},b_{2,m}),&&\sum_{m=1}^\infty (g_2(b_{2,m})-g_1(a_{2,m}))\le \mu_{g_2}^*(E)+\frac{\varepsilon}{2}.
	\end{align*}
	Define $\mathcal R=\{[a_{1,n},b_{1,n})\cap[a_{2,m},b_{2,m}): n,m=1,2,\dots\}\setminus\{\emptyset\}$.
	Observe that the elements of $\mathcal R$ are of the form $[c,d)$, $c,d\in\mathbb R$, $c<d$, since, by construction, we removed those intersections that might be empty. Specifically, denoting $a_n^m=\max\{a_{1,n},a_{2,m}\}$,\ $b_n^m=\min\{b_{1,n},b_{2,m}\}$, $n,m=1,2,\dots$,
	\[\mathcal R=\{[a_n^m,b_n^m)\}_{(n,m)\in\mathcal I},\quad \mathcal I=\left\{(n,m)\in\mathbb N\times\mathbb N: a_n^m<b_n^m\right\}.\]
	Observe that $\mathcal R$ is a countable set. Furthermore, given $x\in E$, since $\mathcal R_1$ and $\mathcal R_2$ are covers of $E$, there exist $n_0,m_0\in\mathbb N$ such that $x\in[a_{1,n_0},b_{1,n_0})$ and $x\in[a_{2,m_0},b_{2,m_0})$, which ensures that $x\in[a_{n_0}^{m_0},b_{n_0}^{m_0})$. This guarantees that $\mathcal R$ is a countable cover of $E$. Therefore,
	\begin{equation}\label{ineqaux1}
	\mu_g^*(E)\le \sum_{(n,m)\in\mathcal I}(g(b_n^m)-g(a_n^m))=\sum_{(n,m)\in\mathcal I}(g_1(b_n^m)-g_1(a_n^m))+\sum_{(n,m)\in\mathcal I}(g_2(b_n^m)-g_2(a_n^m)).
	\end{equation}
	Hence, it is enough to show that
	\begin{equation}\label{outermeasineq}
	\sum_{(n,m)\in\mathcal I}(g_i(b_n^m)-g_i(a_n^m))\le \sum_{k=1}^\infty(g_i(b_{i,k})-g_i(a_{i,k})),\quad i=1,2,
	\end{equation}
	to conclude the result since, in that case, \eqref{ineqaux1} yields
	\begin{align*}
	\mu_g^*(E)\le \sum_{n=1}^\infty(g_1(b_{1,n})-g_1(a_{1,n}))+\sum_{m=1}^\infty(g_2(b_{2,m})-g_2(a_{2,m}))
	<\mu_{g_1}^*(E)+\mu_{g_2}^*(E)+\varepsilon,
	\end{align*}
	which ensures that $\mu_g^*(E)\le \mu_{g_1}^*(E)+\mu_{g_2}^*(E)$ as $\varepsilon>0$ was arbitrarily fixed.

	Let us show that \eqref{outermeasineq} holds.
	For each $n\in\mathbb N$, define $\mathcal I_n=\{m\in\mathbb N:(n,m)\in\mathcal I\}$. Note that $\mathcal I=\bigcup_{n=1}^\infty \mathcal I_n$, and so,
	\begin{equation}\label{sumnm}
	\sum_{(n,m)\in\mathcal I}(g_1(b_n^m)-g_1(a_n^m))=\sum_{n=1}^\infty\left(\sum_{m\in\mathcal{I}_n}(g_1(b_n^m)-g_1(a_n^m))\right).
	\end{equation}
	On the other hand, by definition, we have that $\bigcup_{m\in\mathcal I_n}[a_n^m,b_n^m)\subset [a_{1,n},b_{1,n})$ for each $n\in\mathbb N$.
	Furthermore, since $\mathcal R_2$ is pairwise disjoint, it follows that $\{[a_n^m,b_n^m)\}_{m=1}^\infty$ is also pairwise disjoint for each $n\in\mathbb N$. Hence, for each $n\in\mathbb N$,
	\begin{align*}
	\sum_{m\in\mathcal{I}_n}(g_1(b_n^m)-g_1(a_n^m))=\mu_{g_1}\left(\bigcup_{m\in\mathcal{I}_n}[a_n^m,b_n^m)\right)\le \mu_{g_1}([a_{1,n},b_{1,n}))=g_1(b_{1,n})-g_1(a_{1,n}).
	\end{align*}
	Now \eqref{outermeasineq} for  $i=1$ follows from \eqref{sumnm}. The case $i=2$ is analogous and we omit it.
\qed

In particular, Proposition \ref{gsum} ensures that every set which is $g_i$--measurable for all $i\in\{1,2,\dots,n\}$ is also $\widehat g$--measurable and $\mu_{\widehat g}(E)=\mu_{g_1}(E)+\mu_{g_2}(E)+\dots+\mu_{g_n}(E)$ for all $E\in\bigcap_{i=1}^n \mathcal{LS}_{g_i}$.
The same can be said in regards to the measurability of maps. Furthermore, it follows that if a map $f:X\to\mathbb R$ is $g_i$--integrable for all $i\in\{1,2,\dots,n\}$, then $f$ is $\widehat g$--integrable and
\[\int_X f(s)\dif \widehat g(s)=\int_X f(s)\dif g_1(s)+\int_X f(s)\dif g_2(s)+\dots+\int_X f(s)\dif g_n(s).\]
This fact was implicitly used in \cite{LoMa19} in the proof of the Theorem~4.3, but it was never discussed if such property was true. Here, we have shown that this is the case. This information will be fundamental for the work ahead. Furthermore, it is important to note that throughout this paper, the map $\widehat g$ will denote the map defined as in \eqref{ghatlip}.

\section{The Stieltjes derivative and the concept of continuity}\label{sectionderivative}
In this section we gather some information available in \cite{FriLo17,LoRo14,LoMa19} regarding one of the fundamental tools for this paper: the Stieltjes derivative. Furthermore, we also include some information in those papers regarding different concept of continuity there presented which are also required for the study of differential equations with several Stieltjes derivatives. In particular, we revisit those in \cite{LoMa19} as the results there present some limitations. For the aims of this section, as well as the rest of the paper, we shall assume that $\mathbb R^n$ is endowed with the maximum norm, i.e.,
\[\|x\|=\max\{|x_1|, |x_2|,\dots,|x_n|\},\quad x=(x_1,x_2,\dots,x_n)\in\mathbb R^n.\]

We start by introducing the concept of Stieltjes derivator. From now on, we will refer to nondecreasing and left--continuous maps on $\mathbb R$ as \emph{derivators}. Given a derivator $g$, we will denote by $D_g$ the set of all discontinuity points of $g$. Observe that, given that $g$ is nondecreasing, we can write $D_g=\{ t \in \mathbb R  :  \Delta g(t)>0\}$
where $\Delta g(t)=g(t^+)-g(t)$, $t\in\mathbb R$, and $g(t^+)$ denotes the right handside limit of $g$ at $t$.
With this remark, we now have all the information required to introduce the following definition in \cite{LoRo14}.
\begin{definition}\label{Stieltjesderivative}
	Let $g:\mathbb R\to\mathbb R$ be derivator and $f:\mathbb R\to\mathbb R$. We define the \emph{Stieltjes derivative}, or \emph{$g$--derivative}, of $f$ at a point $t\in \mathbb R\setminus C_g$ as
	\[
	f'_g(t)=\left\{
	\begin{array}{ll}
	\displaystyle \lim_{s \to t}\frac{f(s)-f(t)}{g(s)-g(t)},\quad & t\not\in D_{g},\vspace{0.1cm}\\
	\displaystyle\lim_{s\to t^+}\frac{f(s)-f(t)}{g(s)-g(t)},\quad & t\in D_{g},
	\end{array}
	\right.
	\]
	provided the corresponding limits exist. In that case, we say that $f$ is \emph{$g$--differentiable at $t$}.
\end{definition}
\begin{remark}\label{condgsalto}
	Observe that the points of $C_g$ are excluded from the definition of $g$--derivative. This is because the corresponding limit does not make sense in any neighborhood of these points. Nevertheless, as mentioned before, $\mu_g(C_g)=0$ so, in most cases, this has no impact. Furthermore, note that for $t\in D_g$, $f'_g(t)$ exists if and only if $f(t^+)$ exists and, in that case,
	\[f'_g(t)=\frac{f(t^+)-f(t)}{\Delta g(t)}.\]
\end{remark}

For more information on the Stieltjes derivatives, we refer the reader to \cite{LoRo14,FriLo17}. Here, we will restrict ourselves to the information strictly necessary for the contents of this paper. Along these lines, we include a reformulation of \cite[Theorem~5.4]{LoRo14} with \cite[Definition~~5.1]{LoRo14} added to its statement.
\begin{theorem}\label{gFTC2}
	Let $g:\mathbb R\to\mathbb R$ be derivator
	and $F:[a,b]\to\mathbb R$. The following conditions are equivalent:
	\begin{enumerate}
		\item[\textup{1.}] The function $F$ is $g$--absolutely continuous on $[a,b]$, which we write as $F\in\mathcal{AC}_g([a,b],\mathbb R)$, according to the following definition: for every $\varepsilon>0$, there exists $\delta>0$ such that
		for every open pairwise disjoint family of subintervals $\{(a_n,b_n)\}_{n=1}^m$,
		\begin{displaymath}
		\sum_{n=1}^m (g(b_n)-g(a_n))<\delta\implies
		\sum_{n=1}^m |F(b_n)-F(a_n)|<\varepsilon.
		\end{displaymath}
		\item[\textup{2.}] The function $F$ satisfies the following conditions:
		\begin{enumerate}
			\item[\textup{(i)}] there exists $F'_g(t)$ for $g$--a.a. $t\in[a,b)$;
			\item[\textup{(ii)}] $F'_g\in \mathcal L^1_{g}([a,b),\mathbb R)$, the set of Lebesgue--Stieltjes integrable functions with respect to $\mu_g$;
			\item[\textup{(iii)}] for each $t\in[a,b]$,
			\[F(t)=F(a)+\int_{[a,t)}F'_g(s)\dif g(s).\]
		\end{enumerate}
	\end{enumerate}
\end{theorem}
\begin{remark}\label{gFTC}
	A particularly interesting case of $g$--absolutely continuous function can be found in \textup{Theorem~2.4} and \textup{Proposition 5.2} in \textup{\cite{LoRo14}}. Given $f\in\mathcal L_{g}^1([a,b),\mathbb R)$, the map $F:[a,b]\to\mathbb R$ defined as
	\[F(t)=\int_{[a,t)} f(s)\dif g(s),\]
	is well--defined, $F\in\mathcal{AC}_g([a,b],\mathbb R)$ and $F'_g(t)=f(t)$ for $g$--a.a. $t\in[a,b)$.
\end{remark}

In the work ahead, we will consider systems of differential equations in $\mathbb R^n$ where each component is differentiated with respect to a different derivator. Specifically, we will consider $\boldsymbol g:\mathbb R\to\mathbb R^n$, $\boldsymbol g=(g_1,g_2,\dots, g_n)$, such that each $g_i$, $i\in\{1,2,\dots,n\}$, is a derivator, and we will be looking for solutions on the following set:
\[\mathcal{AC}_{\boldsymbol g}([a,b],\mathbb R^n)= \prod_{i=1}^n \mathcal{AC}_{g_i}([a,b],\mathbb R).\]
That is, we will look for $\boldsymbol g$--absolutely continuous functions (see \cite[Definition 3.4]{LoMa19}), or in other words, functions such that the $i$--th component is $g_i$--absolutely continuous, $i=1,2,\dots,n$.
\begin{remark}\label{pmbgabstohatgabs}
	Observe that for the particular case in which $\boldsymbol g=(g,g,\dots, g)$ for some derivator $g$, we have that $\mathcal{AC}_{\boldsymbol g}([a,b],\mathbb R^n)=\mathcal{AC}_{g}([a,b],\mathbb R^n)$ in the sense presented in \textup{\cite{FriLo17}}.
	Note that, in particular, if $F\in \mathcal{AC}_{\boldsymbol g}([a,b],\mathbb R^n)$, then $F\in \mathcal{AC}_{\widehat g}([a,b],\mathbb R^n)$. To see that, it is enough to note that for each $i\in\{1,2,\dots,n\}$ and each open pairwise disjoint family of subintervals of $[a,b]$, $\{(a_k,b_k)\}_{k=1}^m$,
	\[\sum_{k=1}^m (g_i(b_k)-g_i(a_k))\le \sum_{k=1}^m (\widehat g(b_k)-\widehat g(a_k)).\]
\end{remark}

Theorem \ref{gFTC2} ensures that, given $F\in \mathcal{AC}_{\boldsymbol g}([a,b],\mathbb R^n)$, $F=(F_1,F_2,\dots,F_n)$,  $(F_i)'_{g_i}(t)$ exists for $g_i$--a.a. $t\in[a,b)$ and $(F_i)'_{g_i}\in\mathcal L^1_{g_i}([a,b),\mathbb R)$, $i\in\{1,2,\dots,n\}$. Moreover, we have that
\[F_i(t)=F_i(a)+\int_{[a,t)} (F_i)'_{g_i}(s)\dif g_i(s),\quad t\in[a,b],\quad  i=1,2,\dots, n.\]
Throughout this paper, we will use ``component-by-component'' notation for the derivatives and integrals so that expressions such as the previous one can simply be reduced to
\[F(t)=F(a)+\int_{[a,t)} F'_{\boldsymbol g}(s) \dif \boldsymbol g(s),\quad t\in[a,b].\]

The rest of this section is dedicated to the study of the concept of continuity with respect to a map  $\boldsymbol g:\mathbb R\to\mathbb R^n$, $\boldsymbol g=(g_1,g_2,\dots, g_n)$, such that each $g_i$, $i\in\{1,2,\dots,n\}$, is a derivator. Some research on this topic was carried out in \cite{FriLo17} in the setting of a unique derivator, and in \cite{LoMa19} in the general setting. Unfortunately, some of the results in the latter are incorrect. Here, we show their limitations and amend some of those errors. We start by introducing the following concept that contains \cite[Definition 3.1]{FriLo17}.
\begin{definition}\label{multicont}
	Let $g:\mathbb R\to\mathbb R$ and $\boldsymbol g:\mathbb R\to\mathbb R^n$, $\boldsymbol g=(g_1,g_2,\dots, g_n)$, be such that $g$, $g_i$, $i=1,2,\dots,n$, are derivators; and consider $f:A\subset \mathbb R\to \mathbb R^n$, $f=(f_1,f_2,\dots, f_n)$. We say that $f$ is \emph{$g$--continuous at $t\in A$} if for every $\varepsilon>0$, there exists $\delta>0$ such that
		\[\|f(s)-f(t)\|<\varepsilon,\quad\mbox{for all }s\in A\mbox{ such that }\ |g(s)-g(t)|<\delta.\]
	We say that $f$ is \emph{$\boldsymbol g$--continuous at $t\in A$} if $f_i$ is $g_i$--continuous at $t$ for each $i\in\{1,2,\dots,n\}$.
	If $f$ is $\boldsymbol g$--continuous at every $t \in A$, we say that $f$  is \emph{$\boldsymbol g$--continuous on $A$}.
\end{definition}
\begin{remark}\label{gabsbound}
	We denote by $\mathcal C_{\boldsymbol g}([a,b],\mathbb R^n)$ the set of $\boldsymbol g$--continuous functions on $[a,b]$ with values in $\mathbb R^n$ and by $\mathcal{BC}_{\boldsymbol g}([a,b],\mathbb R^n)$ the set of functions in $\mathcal C_{\boldsymbol g}([a,b],\mathbb R^n)$ which are also bounded. Note that, by definition, we have that
	\[\mathcal{BC}_{\boldsymbol g}([a,b],\mathbb R^n)= \prod_{i=1}^n \mathcal{BC}_{g_i}([a,b],\mathbb R),\]
	which, together with \textup{\cite[Proposition 5.5]{FriLo17}}, ensures that $\mathcal{AC}_{\boldsymbol g}([a,b],\mathbb R^n)\subset\mathcal{BC}_{\boldsymbol g}([a,b],\mathbb R^n)$.
\end{remark}

The following result, which can be directly deduced from \cite[Proposition 3.2]{FriLo17}, contains some basic properties regarding $\boldsymbol g$--continuous functions that will be useful in what follows.
\begin{proposition}\label{multigcontprop}
	Let $f:[a,b]\to \mathbb R^n$, $f=(f_1,f_2,\dots, f_n)$,  be a $\boldsymbol g$--continuous map on $[a,b]$. Then:
	\begin{enumerate}
		\item[\textup{1.}] For each $i\in\{1,2,\dots,n\}$, $f_i$ is continuous from the left at every $t \in (a,b]$.
		\item[\textup{2.}] For each $i\in\{1,2,\dots,n\}$, $f_i$ is continuous at every $t \in (a,b)$ at which $g_i$ is continuous.
		\item[\textup{3.}] For each $i\in\{1,2,\dots,n\}$, $f_i$ is constant on every $[c,d] \subset [a,b]$ at which $g_i$ is constant.
	\end{enumerate}
\end{proposition}

Proposition~\ref{multigcontprop} allows us to establish some relations between the continuity of a derivator with respect to another one and some of their characteristic sets, as presented in the next result.
\begin{proposition}\label{CgiDgi}
	Let $g_1,g_2:\mathbb R\to\mathbb R$ be derivators. Then $g_1$ is $g_2$--continuous if and only if $C_{g_2}\subset C_{g_1}$ and $D_{g_1}\subset D_{g_2}$.
\end{proposition}
{\bf Proof.}
First, assume that $g_1$ is $g_2$--continuous. If $t\in C_{g_2}$, then there exists $\delta>0$ such that $g_2$ is constant in $(t-\delta,t+\delta)$. In particular, $g_2$ is constant on $[t-\delta/2,t+\delta/2]$. Proposition~\ref{multigcontprop} guarantees that $g_1$ is constant on that same interval, so it is constant on $(t-\varepsilon,t+\varepsilon)$, $\varepsilon\in(0,\delta/2)$, i.e., $t\in C_{g_1}$.

On the other hand, if $t\in\mathbb R\setminus D_{g_2}$, we have that $g_2$ is continuous at $t$. In that case,
Proposition~\ref{multigcontprop} yields that $g_1$ is continuous at $t$, and so $t\in \mathbb R\setminus D_{g_1}$. Thus, $D_{g_1}\subset D_{g_2}$.

Conversely, assume that $C_{g_2}\subset C_{g_1}$ and $D_{g_1}\subset D_{g_2}$ and let $t\in\mathbb R$. Observe that showing that $g_1$ is $g_2$--continuous at $t$ is equivalent to showing that for each $\varepsilon>0$, there exist $\delta_1,\delta_2>0$ such that
\begin{align}
&0\le g_1(s)-g_1(t)<\varepsilon,&\mbox{for all }s\ge t\mbox{ such that }0\le g_2(s)-g_2(t)<\delta_1,\label{rightcont}\\
&0\le g_1(t)-g_1(s)<\varepsilon,&\mbox{for all }s<t\mbox{ such that }0\le g_2(t)-g_2(s)<\delta_2\label{leftcont}.
\end{align}

Let $\varepsilon>0$. In order to show that \eqref{rightcont} holds, define $A_t=\{s\in[t,+\infty): g_2(s)=g_2(t)\}\not=\emptyset$. If $A_t$ is not bounded from above, we have that $g_2(s)=g_2(t)$ for all $s\in[t,+\infty)$.
In that case, $(t,+\infty)\subset C_{g_2}\subset C_{g_1}$, which means that $g_1(s)=g_1(t)$ for all $s\in(t,+\infty)$, and so \eqref{rightcont} holds for any $\delta_1>0$.
Otherwise, $A_t$ is bounded from above and we can define $a_t=\sup A_t\in[t,+\infty)$. Observe that, by definition, we have that $g_2(s)=g_2(t)$ for all $s\in[t,a_t)$ and $g_2(s)>g_2(t)$ for all $s\in(a_t,+\infty)$. In particular, we have that $(t,a_t)\subset C_{g_2}\subset C_{g_1}$, which implies that $g_1(s)=g_1(t)$ for all $s\in[t,a_t)$. Now, since $g_1$ and $g_2$ are left--continuous, it follows that $g_1(s)=g_1(t)$ and $g_2(s)=g_2(t)$ for all $s\in[t,a_t]$.

Now, if $a_t\in D_{g_1}$, we have that $a_t\in D_{g_2}$. In that case, taking $\delta_1=\Delta g_2(a_t)>0$, we have that $\{s\ge t:0\le g_2(s)-g_2(t)<\delta_1 \}=[t,a_t]$, and so \eqref{rightcont} follows (even when $a_t=t$, as it becomes trivial). Otherwise, we have that $a_t\not\in D_{g_1}$ which means that $g_1$ is continuous from the right at $a_t$, so there exists $\delta'>0$ such that
\begin{equation}\label{g1rightcont}
0\le g_1(s)-g_1(a_t)<\varepsilon\quad\mbox{for all }s\ge a_t\mbox{ such that }0\le s-a_t<\delta'.
\end{equation}
Observe that it is enough to show that there exists $\delta_1>0$ such that
\begin{equation}\label{claim1}
0\le s-a_t<\delta'\quad \mbox{for all }s\ge a_t\mbox{ such that }0\le g_2(s)-g_2(a_t)<\delta_1,
\end{equation}
to obtain that \eqref{rightcont} holds. Indeed, if $a_t=t$, then \eqref{rightcont} follows directly from \eqref{g1rightcont} and \eqref{claim1}. Otherwise, $a_t>t$ and, in that case, if $s\ge t$ is such that $0\le g_2(s)-g_2(t)<\delta_1$ we have that either $s\in [t,a_t]$, implying that $g_1(s)-g_1(t)=0$, or $s>a_t$ and, in that case, it follows again from \eqref{g1rightcont} and \eqref{claim1}.

Reasoning by contradiction, if \eqref{claim1} does not hold, then for every $n\in{\mathbb N}$ there exists $s_n\ge a_t$ such that $s_n-a_t\ge \delta'$ and $0\le g_2(s_n)-g_2(a_t)<1/n$. Let $s=\inf\{s_n\}_{n\in{\mathbb N}}\ge a_t+\delta'>a_t$. Since $g_2$ is nondecreasing, $0\le g_2(s)-g_2(a_t)\le g_2(s_n)-g_2(a_t)<1/n$ for every ${n\in{\mathbb N}}$. Therefore, $g_2(s)=g_2(a_t)$ which contradicts the fact that $g_2(s)>g_2(t)=g_2(a_t)$ for all $s\in[a_t,+\infty)$.

Finally, for \eqref{leftcont} we proceed in a similar manner. Define $B_t=\{s\in(-\infty,t]: g_2(s)=g_2(t)\}\not=\emptyset$. If $B_t$ is not bounded from below, an analogous reasoning to the one for $A_t$ shows that \eqref{leftcont} holds trivially. Otherwise, $B_t$ is bounded from below and we can define
$b_t=\inf B_t\in(-\infty, t]$. Hence, by definition, we have that $g_2(s)=g_2(t)$ for all $s\in(b_t,t]$ and $g_2(s)<g_2(t)$ for all $s\in(-\infty, b_t]$.  In particular, $(b_t,t)\subset C_{g_2}\subset C_{g_1}$, which guarantees that $g_1(s)=g_1(t)$ for all $s\in(b_t,t]$.

If $b_t\in D_{g_2}$, then taking $\delta_2=\Delta g_2(b_t)>0$, we have that $\{s<t:0\le g_2(s)-g_2(t)<\delta_1 \}=(b_t,t)$, and so \eqref{rightcont} follows (even when $b_t=t$, as it becomes vacuous). Otherwise, $b_t\not\in D_{g_2}$ which implies that $b_t\not\in D_{g_1}$. Observe that this ensures that
\begin{equation}\label{condclaim2}
g_1(s)=g_1(t)\ \mbox{ and }\ g_2(s)=g_2(t),\quad\mbox{for all $s\in[b_t,t]$.}
\end{equation}
Furthermore, since $g_1$ is left--continuous at $b_t$, there exists $\delta''>0$ such that
\begin{equation*}
0\le g_1(b_t)-g_1(s)<\varepsilon\quad\mbox{for all }s<b_t\mbox{ such that }0<b_t-s<\delta''.
\end{equation*}
Now, an analogous reasoning to the one for \eqref{claim1} shows that, given that \eqref{condclaim2} holds, there exists $\delta_2>0$ such that $0<b_t-s<\delta''$ for all $s<b_t$ satisfying that $0<g_2(b_t)-g_2(s)<\delta_2$, which is enough to conclude that \eqref{leftcont} holds in a similar fashion to \eqref{rightcont}. \qed
\begin{remark}
	Condition \eqref{rightcont} can be interpreted as ``being $g_2$--continuous at $t$ from the right'' and, similarly, \eqref{leftcont} as ``being $g_2$--continuous at $t$ from the left''.
\end{remark}

As pointed out in \cite{FriLo17}, the concept of $g$--continuity can be understood as continuity in the topological sense. Specifically, a map $f:\mathbb R\to\mathbb R$ is $g$--continuous if $f:(\mathbb R,\tau_g)\to(\mathbb R,\tau_u)$ is continuous, where $\tau_u$ is the usual topology of $\mathbb R$ and $\tau_g$ is the topology generated by the sets
\begin{equation}\label{gopenball}
B_{g}(x,r):=\{y\in \mathbb R\ :\ |g(y)-g(x)|<r\},\quad x\in\mathbb R,\ r>0.
\end{equation}
With this charaterization, observe that $g_1$ is $g_2$--continuous if and only if $\tau_{g_1}\subset\tau_{g_2}$.
Indeed, first observe that $g_1$ is trivially $g_1$--continuous. This means that $g_1^{-1}(U)\in\tau_{g_1}$ for any $U\in \tau_u$. Hence, if $\tau_{g_1}\subset\tau_{g_2}$, we have that $g_1^{-1}(U)\in\tau_{g_2}$ for any $U\in \tau_u$, or equivalently, $g_1$ is $g_2$--continuous. Conversely, if $g_1$ is $g_2$--continuous, and given that the corresponding topologies are generated by the sets in \eqref{gopenball}, it is enough to show that for any $t\in\mathbb R$, $r>0$ and $s\in B_{g_1}(t,r)$, there exists $r_s>0$ such that $B_{g_2}(s,r_s)\subset B_{g_1}(t,r)$. Given $t\in\mathbb R$, $r>0$ and $s\in B_{g_1}(t,r)$, define $\delta=r-|g_1(t)-g_1(s)|>0$.  Then, $B_{g_1}(s,\delta)\subset B_{g_1}(t,r)$ and, since $g_1$ is $g_2$--continuous, there exists $r_s>0$ such that
\[|g_1(z)-g_1(s)|<\delta,\quad\mbox{for all }z\in\mathbb R\mbox{ such that }|g_2(z)-g_2(s)|<r_s.\]
Hence, if $z\in B_{g_2}(s,r_s)$, we have that $z\in B_{g_1}(s,\delta)\subset B_{g_1}(t,r)$, which finishes the proof of the equivalence. Combining this with Proposition \ref{CgiDgi}, we obtain the following result.
\begin{corollary}\label{corit}
Let $g_1,g_2:\mathbb R\to\mathbb R$ be derivators. Then the following statements are equivalent:
	\begin{enumerate}
		\item[\textup{1.}] $g_1$ is $g_2$--continuous and $g_2$ is $g_1$--continuous.
		\item[\textup{2.}] $C_{g_2}= C_{g_1}$ and $D_{g_1}= D_{g_2}$.
		\item[\textup{3.}] $\tau_{g_1}=\tau_{g_2}$.
	\end{enumerate}
\end{corollary}

As we will see later, this result will allow us to relate the concept of $\boldsymbol g$--continuity with a similar type of continuity in \cite{LoMa19}. Furthermore, Corollary~\ref{corit} implies that the topologies of derivators can be classified in terms of the sets $C_g$ and $D_g$ as the following theorem shows.

\begin{theorem}[Classification of derivator topologies] \label{thmcdp}
	Let $G$ be the family of all derivators and define the following equivalence relation on $G$:
	\[g_1\sim g_2:\iff \tau_{g_1}=\tau_{g_2}.\]
	Let ${\mathcal G}:=G|_\sim$ and $\mathcal H=\{(C,D)\in\mathcal P(\mathbb R)\times \mathcal P(\mathbb R): C\mbox{ open, } D\mbox{ at most countable, }C\cap D=\emptyset\}$. Then, the map $\Phi:{\mathcal G}\to{\mathcal H}$, $\Phi([g])=(C_g,D_g)$, is a bijection.
\end{theorem}
{\bf Proof.}
	First, observe Corollary~\ref{corit} ensures that $\Phi$ is well defined and injective. To see that it is surjective, let us denote by $m$ the Lebesgue measure on $\mathbb R$ and $\llbracket a,b\rrbracket:=[\min\{a,b\},\max\{a,b\}]$, $a,b\in\mathbb R$. Then, given $(C,D)\in{\mathcal H}$, let $D=\{d_n\}_{n\in\Lambda}$, $\Lambda\subset{\mathbb N}$, and define $g_{_{C,D}}:{\mathbb R}\to{\mathbb R}$ as
    \[ g_{_{C,D}}(t):=\operatorname{sgn}(t) m\left(\llbracket 0,t\rrbracket{\backslash} C\right)+\sum_{\substack{d_n\in D\\ d_n<t}}2^{-n}.\]
    Note that $\llbracket 0,t\rrbracket{\backslash} C$ is a closed set, thus Lebesgue--measurable. Hence, $g_{_{C,D}}$ is well defined. Furthermore, $g_{_{C,D}}$ is left--continuous and nondecreasing, $D_g=D$ and $C_g=C$, which shows that $\Phi$ is surjective.\qed

\begin{example}
	To illustrate \textup{Theorem~\ref{thmcdp}} we consider $[a,b]=[0,1]$ and the class of the derivator \[g(t)=t+\sum_{\substack{\frac{1}{n}<t\\ n\in\mathbb N}}2^{-n}.\]
	The class $[g]$ is the set of all derivators that generate the same topology. Observe that $C_g=\emptyset$ and $D_g=\{\frac{1}{n}\}_{n\in\mathbb N}$, so $[g]$ is characterized by $\left(\emptyset,\{\frac{1}{n}\}_{n\in\mathbb N}\right)\in\mathcal H$. Thus, $[g]$ can be explicitly expressed as	\[[g]=\left\{h(t)=\phi(t)+\sum_{\substack{\frac{1}{n}<t\\ n\in\mathbb N}}\alpha_n\ :\ \alpha_n>0,\ \sum_{n\in\mathbb N}\alpha_n<\infty,\ \phi \text{ strictly increasing and continuous}\right\}.\]
	\end{example}
As mentioned before, in \cite[Definition 3.1]{LoMa19} we find a definition of continuity with respect to $\boldsymbol g$ alternative to the one provided in Definition~\ref{multicont}. For completeness, we include \cite[Definition 3.1]{LoMa19} before comparing the two concepts.

\begin{definition}\label{vecmulticont}
	Let $\boldsymbol g:\mathbb R\to\mathbb R^n$, $\boldsymbol g=(g_1,g_2,\dots, g_n)$, be such that each $g_i$, $i\in\{1,2,\dots,n\}$, is a derivator. A function $f:A\subset \mathbb R\to \mathbb R^n$, $f=(f_1,f_2,\dots, f_n)$, is \emph{$\vec{\boldsymbol{g}}$--continuous at point $t\in A$} if, for every $\varepsilon>0$, there exists $\delta>0$ such that
	\[\|  f(t)-  f(s)\|<\varepsilon,\quad \mbox{for all }s\in A\mbox{ such that }\ \|\boldsymbol g(t)-\boldsymbol g(s)\|<\delta.\]
	If it is $\vec{\boldsymbol{g}}$--continuous at every point $t \in A$, we say that $f$  is \emph{$\vec{\boldsymbol{g}}$--continuous on $A$}.
\end{definition}

In \cite{LoMa19}, the authors claimed that Definitions \ref{multicont} and \ref{vecmulticont} are equivalent. Nevertheless, a careful reader might notice that the proof only shows that $\boldsymbol g$--continuity implies $\vec{\boldsymbol{g}}$--continuity. Furthermore, the reverse implication is not true, as shown in the next example.

\begin{example}\label{exvecnomulti}
	Consider $\boldsymbol g, f:\mathbb R\to\mathbb R^2$, $\boldsymbol g=(g_1,g_2)$, $f=(f_1,f_2)$, defined as
	\begin{equation}\label{eqvecnomulti}
	\boldsymbol g(t)=\left\{
	\begin{array}{ll}
	(0,t),&\mbox{if }t\le 0,\vspace{0.1cm}\\
	\displaystyle\left(0,t+1\right),& \mbox{if }t>0,
	\end{array}
	\right.
	\quad \quad \quad
	f(t)=\left\{
	\begin{array}{ll}
	(t,0),&\mbox{if }t\le 0,\vspace{0.1cm}\\
	\displaystyle\left(t+1,\frac{\sin(1/t)}{t}\right),& \mbox{if }t>0.
	\end{array}
	\right.
	\end{equation}
	Note that $f$ cannot be $\boldsymbol g$--continuous as $f_1$ is not constant, see \textup{Proposition \ref{multigcontprop}}. However, $f$ is $\vec{\boldsymbol{g}}$--continuous. Indeed, first note that $\|\boldsymbol g(t)-\boldsymbol g(s)\|=|g_2(t)-g_2(s)|$, $s,t\in \mathbb R$.
	Thus, showing that $f$ is $\vec{\boldsymbol{g}}$--continuous is equivalent to showing that $f$ is $g_2$--continuous. Since we are considering the $\max$--norm in $\mathbb R^n$, it suffices to show that $f_1$ and $f_2$ are $g_2$--continuous. Now, $f_1$ is trivially $g_2$--continuous as $f_1=g_2$, and \textup{\cite[Example~3.3]{FriLo17}} shows that $f_2$ is $g_2$--continuous.
\end{example}

It is important to note that the misinformation about the relations between Definitions \ref{multicont} and \ref{vecmulticont} does not affect the existence and uniqueness results in \cite{LoMa19}. However, it has some consequences when it comes to the study of solutions in the classical sense. Specifically, this affects Proposition 4.6 and Theorem 4.8 in \cite{LoMa19}. In the next section we discuss the implications of these facts. Nevertheless, in some contexts, those results remain true as a consequence of the following result.
\begin{proposition}\label{equivpmbgvecg}
	Let $\boldsymbol g:\mathbb R\to\mathbb R^n$, $\boldsymbol g=(g_1,g_2,\dots, g_n)$, be such that each $g_i$, $i\in\{1,2,\dots,n\}$, is a derivator. Then, the following are equivalent:
	\begin{enumerate}
		\item[\textup{(i)}] Every $\vec{\boldsymbol{g}}$--continuous map is $\boldsymbol g$--continuous.
		\item[\textup{(ii)}] For each $j,k\in\{1,2,\dots,n\}$, the map $g_k:\mathbb R\to\mathbb R$ is $g_j$--continuous.
		\item[\textup{(iii)}] For each $j,k\in\{1,2,\dots,n\}$, $C_{g_j}= C_{g_k}$ and $D_{g_j}= D_{g_k}$.
		\item[\textup{(iv)}] For each $j,k\in\{1,2,\dots,n\}$, $\tau_{g_j}=\tau_{g_k}$.
	\end{enumerate}
\end{proposition}
\noindent
{\bf Proof.}
Assume that (i) holds. Fix $k\in\{1,2,\dots,n\}$ and consider the map \[G(t)=(g_k(t),g_k(t),\dots,g_k(t)),\quad t\in\mathbb R.\]
Observe that it is enough to show that $G$ is $\vec{\boldsymbol{g}}$--continuous to prove (ii). This is straightforward. Indeed, let $t\in\mathbb R$, $\varepsilon>0$ and take $\delta=\varepsilon$. Then, if $s\in\mathbb R$ is such that $\|\boldsymbol g(t)-\boldsymbol g(s)\|<\delta$, it follows that
\[\|G(t)-G(s)\|=|g_k(t)-g_k(s)|\le\|\boldsymbol g(t)-\boldsymbol g(s)\|<\delta=\varepsilon.\]

Conversely, assume that (ii) holds and let $f:A\subset \mathbb R\to\mathbb R^n$, $f=(f_1,f_2,\dots,f_n)$, be a $\vec{\boldsymbol{g}}$--continuous map. Fix $j\in\{1,2,\dots,n\}$, $t\in A$ and $\varepsilon>0$. Since $f$ is $\vec{\boldsymbol{g}}$--continuous, there exists $\gamma>0$ such that
\[\|f(t)-f(s)\|<\varepsilon,\quad\mbox{for all $s\in A$ such that }\|\boldsymbol g(t)-\boldsymbol g(s)\|<\gamma.\]
On the other hand, for each $k\in\{1,2,\dots, n\}$, $g_k$ is $g_j$--continuous and so, there exists $\delta_k>0$ such that
\[|g_k(t)-g_k(s)|<\gamma,\quad \mbox{for all $s\in \mathbb R$ such that } |g_j(t)-g_j(s)|<\delta_k.\]
Take $\delta=\min\{\delta_1,\delta_2,\dots,\delta_n\}$. In that case, if $s\in A$ such that $|g_j(t)-g_j(s)|<\delta$, we have that $\|\boldsymbol g(t)-\boldsymbol g(s)\|<\gamma$, which ensures that $|f_j(t)-f_j(s)|\le \|f(t)-f(s)\|<\varepsilon$. Hence, $f_j$ is $g_j$--continuous and, since $j\in\{1,2,\dots,n\}$ was arbitrarily fixed, $f$ is $\boldsymbol g$--continuous.

The rest of the result is a consequence of Corollary~\ref{corit}.
\qed

Interestingly enough, it is possible to establish some other connections between the continuity in the sense of Definition \ref{vecmulticont} and \cite[Definition~3.1]{FriLo17} for an adequate choice of a derivator.

\begin{proposition}\label{multiconttocont}
	Let $\boldsymbol g:\mathbb R\to\mathbb R^n$, $\boldsymbol g=(g_1,g_2,\dots, g_n)$, be such that each $g_i$, $i\in\{1,2,\dots,n\}$, is a derivator.
	Then, $f:A\subset \mathbb R\to \mathbb R^n$ is $\vec{\boldsymbol{g}}$--continuous at $t\in A$ if and only if $f$ is $\widehat g$--continuous at $t$.
\end{proposition}
\noindent
{\bf Proof.}
Observe that for all $s\in A$,
\begin{equation}\label{vecghatgineq}
\|\boldsymbol g(t)-\boldsymbol g(s)\|\le |\widehat g(t)-\widehat g(s)|\le n\|\boldsymbol g(t)-\boldsymbol g(s)\|.
\end{equation}
Indeed, it follows from the triangular inequality that $|\widehat g(t)-\widehat g(s)|\le n\|\boldsymbol g(t)-\boldsymbol g(s)\|$ for any $s\in A$.
For the other inequality we consider two cases. If $t\ge s$, since $\widehat g$ and each $g_i$ are nondecreasing, we have
\begin{equation}\label{eqaux1}
|\widehat g(t)-\widehat g(s)|=\widehat g(t)-\widehat g(s)=\sum_{i=1}^n (g_i(t)-g_i(s))\ge \max_{i=1,\dots,n}\{g_i(t)-g_i(s)\}=\|\boldsymbol g(t)-\boldsymbol g(s)\|.
\end{equation}
On the other hand, if $t<s$, we proceed analogously to \eqref{eqaux1}, interchanging the roles of $t$ and $s$.
Hence, \eqref{vecghatgineq} holds. Now, the equivalence between the two types of continuity follows.
\qed

Proposition \ref{multiconttocont} not only provides a simple condition to check if a map is $\vec{\boldsymbol{g}}$--continuous, but we can also deduce some interesting properties for this type of maps through the results in \cite{FriLo17}. In particular, \cite[Corollary 3.5]{FriLo17} yields the following result that is fundamental for the uniqueness results in the following section.
\begin{proposition}\label{multiborel}
	Let $A\subset\mathbb R$ be a Borel set and $f:A\subset \mathbb R\to\mathbb R^n$
	be $\vec{\boldsymbol{g}}$--continuous on $A$. Then, $f$ is Borel measurable.
\end{proposition}
\begin{remark}\label{multigmeasurable}
	Since every Borel set is Lebesgue--Stieltjes measurable, it follows that every Borel measurable map is Lebesgue--Stieltjes measurable. In particular, we have that if $A$ is a Borel set and $f:A\subset \mathbb R\to\mathbb R^n$ is $\vec{\boldsymbol{g}}$--continuous on $A$, then $f$ is $g_j$--measurable for all $j\in\{1,2,\dots,n\}$.
\end{remark}

Another type of continuity defined in terms of $\boldsymbol g$ that was introduced in \cite{LoMa19} is what the authors called $(\vec{\boldsymbol{g}}\times \Id)$--continuity, which is a generalization of \cite[Definition 7.7]{FriLo17}. These concepts of continuity were introduced for the study of classical solutions in both papers. As mentioned earlier, the results in \cite{LoMa19} are partially incorrect and we aim to correct them in this paper. To that end, we introduce the following definition, which is an alternative generalization of \cite[Definition~7.7]{FriLo17}.

\begin{definition}\label{multigIdcont}
	Let $g:\mathbb R\to\mathbb R$ and $\boldsymbol g:\mathbb R\to\mathbb R^n$, $\boldsymbol g=(g_1,g_2,\dots, g_n)$, be such that $g$, $g_i$, $i=1,2,\dots,n$, are derivators; and consider $f: A\times B \subset \mathbb R\times \mathbb R^{n} \to \mathbb R^n$, $f=(f_1,f_2,\dots,f_n)$.
	We say that $f$ is \emph{$(g \times \Id)$--continuous at $(t,x)\in A\times B$} if for every $\varepsilon>0$, there exists $\delta>0$ such that
	\[\|f(s,y)-f(t,x)\|<\varepsilon\quad\mbox{for all }(s,y)\in A\times B\mbox{ such that }\ |g(s)-g(t)|<\delta\ \mbox{ and }\ \|y-x\|<\delta.\]
	We say that $f$ is \emph{$(\boldsymbol g \times \Id)$--continuous at $(t,x)$} if each $f_i$ is $(g_i\times \Id)$--continuous at $(t,x)$, $i\in\{1,2,\dots,n\}$. We  say that $f$ is \emph{$(\boldsymbol g \times\Id)$--continuous in $A\times B$} if it is $(\boldsymbol g \times\Id)$--continuous at every $(t,x) \in A\times B$.
\end{definition}
On the other hand, the corresponding definition in \cite{LoMa19} reads as follows.
\begin{definition}\label{vecgIdcont}
	Let $\boldsymbol g:\mathbb R\to\mathbb R^n$, $\boldsymbol g=(g_1,g_2,\dots, g_n)$, be such that each $g_i$, $i\in\{1,2,\dots,n\}$, is a derivator. A function $f: A \times B \subset \mathbb R\times \mathbb R^{n} \to \mathbb R^n$, $f=(f_1,f_2,\dots,f_n)$ is $(\vec{\boldsymbol{g}} \times \Id)$--continuous at $(t,x)\in A\times B$ if for every $\varepsilon>0$, there exists $\delta>0$ such that
	\[\|f(s,y)-f(t,x)\|<\varepsilon\quad\mbox{for all }(s,y)\in A\times B\mbox{ such that }\ \|\boldsymbol g(s)-\boldsymbol g(t)\|<\delta\ \mbox{ and }\ \|y-x\|<\delta.\]
	We  say that $f$ is $(\vec{\boldsymbol{g}} \times\Id)$--continuous in $A\times B$ if it is $(\vec{\boldsymbol{g}} \times\Id)$--continuous at every $(t,x) \in A\times B$.
\end{definition}

The relations between Definitions \ref{multigIdcont} and \ref{vecgIdcont} are analogous to the ones between Definitions \ref{multicont} and \ref{vecmulticont}. In particular, we have the following result.
\begin{proposition}\label{pmbgIdtovecgId}
	Let $f: A\times B \subset \mathbb R\times \mathbb R^{n} \to \mathbb R^n$, $f=(f_1,f_2,\dots, f_n)$. If $f$ is $(\boldsymbol g\times\Id)$--continuous at $(t,x)\in A\times B$, then $f$ is $(\vec{\boldsymbol{g}}\times\Id)$--continuous at $(t,x)$.
\end{proposition}
\noindent
{\bf Proof.}
	Fix $\varepsilon>0$. Given $i\in\{1,2,\dots,n\}$, we have that $f_i$ is $(g_i\times\Id)$--continuous at $(t,x)$, so there exists $\delta_i>0$ such that
	\[|f_i(s,y)-f_i(t,x)|<\varepsilon\quad\mbox{for all }(s,y)\in A\times B\mbox{ such that }\ |g_i(s)-g_i(t)|<\delta_i\ \mbox{ and }\ \|y-x\|<\delta_i.\]
	Take $\delta=\min\{\delta_1,\delta_2,\dots,\delta_n\}.$ Then, if $\|\boldsymbol g(t)-\boldsymbol g(s)\|<\delta$ and $\|  x-  y\| < \delta_i$, we have that $|g_i(t)-g_i(s)|<\delta_i$ and $\|  x-  y\| < \delta_i$ for all $i\in\{1,2,\dots,n\}$, which implies that $\|f(t,x)-  f(s,y)\|<\varepsilon$.
\qed

Once again, the reverse implication does not hold. To see that this is the case, it is enough to consider $\boldsymbol g$ and $f$ as in \eqref{eqvecnomulti} and $F:\mathbb R\times\mathbb R^2\to\mathbb R^2$ defined as $F(t,(x,y))=f(t)$, and note $F$ does not depend on $(x,y)$, which implies that $(\boldsymbol g\times \Id)$ and $(\vec{\boldsymbol{g}}\times \Id)$--continuity reduce to $\boldsymbol g$ and $\vec{\boldsymbol{g}}$--continuity. Furthermore, using a  similar reasoning, we can obtain analogous results to Propositions \ref{equivpmbgvecg} and \ref{multiconttocont} in the context of Definitions \ref{multigIdcont} and \ref{vecgIdcont}.

Interestingly enough, within the proof of \cite[Theorem 4.8]{LoMa19} the authors proved correctly the following superposition result involving  Definitions \ref{vecmulticont} and \ref{vecgIdcont}.

\begin{lemma}\label{pmbgcomp}
	Let $f: A\times B \subset \mathbb R\times \mathbb R^{n} \to\mathbb R^n$ be a $(\vec{\boldsymbol{g}}\times \Id)$--continuous function on $A\times B$. If $x:A\to B$ is $\vec{\boldsymbol{g}}$--continuous on $A$, then the composition $f(\cdot,x(\cdot))$ is $\vec{\boldsymbol{g}}$--continuous on $A$.
\end{lemma}
\begin{remark}\label{remgcompnonbien}
	\textup{Lemma \ref{pmbgcomp}} guarantees that the composition of a $\boldsymbol g$--continuous map with a $(\boldsymbol g\times\Id)$--continuous one is $\vec{\boldsymbol{g}}$--continuous. Nevertheless, we cannot assure that the composition is $\boldsymbol g$--continuous.
	Indeed, consider $\boldsymbol g:\mathbb R\to\mathbb R^2$, $\boldsymbol g=(g_1,g_2)$, be such that $g_1$, $g_2$ are derivators satisfying $\Delta g_1(t_0)=0$ and $\Delta g_2(t_0)>0$ for some $t_0\in\mathbb R$. Let $I$ be a neighborhood of $t_0$ and consider the maps $x:I\to\mathbb R^2$, $f:I\times \mathbb R^2\to\mathbb R^2$, $f=(f_1,f_2)$, defined as
	\[x(t)=(g_1(t),g_2(t)),\quad\quad\quad f(t,(y,z))=(g_1(t)-z,g_2(t)-y).\]
	It is clear that $x$ is $\boldsymbol g$--continuous at $t_0$. Furthermore, observe that
	\[\lim_{t\to t_0^+} f_1(t,x(t))=\lim_{t\to t_0^+}g_1(t)-g_2(t)=g_1(t_0)-g_2(t_0^+)<g_1(t_0)-g_2(t_0)=f_1(t_0,x(t_0)).\]
	This implies that the map $f(\cdot,x(\cdot))$ is not $\boldsymbol g$--continuous at $t_0$, as $f_1(\cdot,x(\cdot))$ is not $g_1$--continuous at $t_0$, see \textup{Proposition \ref{multigcontprop}}. However, the map $f$ is $(\boldsymbol g\times \Id)$--continuous at $(t_0,(g_1(t_0),g_2(t_0)))$. Indeed, we shall only show that $f_1$ is $(g_1\times \Id)$--continuous at $t_0$ as the proof for $f_2$ being $(g_2\times \Id)$--continuous is analogous.
	Let $\varepsilon>0$ and take $0<\delta<\varepsilon/2$. Denote $u_0=(g_1(t_0),g_2(t_0))$. If $(t,(x,y))\in I\times\mathbb R^2$ is such that $|g_1(t_0)-g_1(t)|<\delta$ and $\|u_0-(y,z)\|<\delta$, then
	\begin{align*}
	|f_1(t_0,u_0)-f_1(t,(y,z))|\le &  |g_1(t_0)-g_1(t)|+|y-g_2(t_0)|\\ \le &  |g_1(t_0)-g_1(t)|+\|u_0-(x,y)\|<2\delta<\varepsilon.
	\end{align*}
\end{remark}

\section{The initial value problem with several derivativors}\label{sectionivp}

We now turn our attention to the study of initial value problems in the context of $\boldsymbol g$--differential equations. That is, given $\boldsymbol g:\mathbb R\to\mathbb R^n$, $\boldsymbol g=(g_1,g_2,\dots, g_n)$, such that each $g_i$, $i\in\{1,2,\dots,n\}$, is a derivator, we will study problems of the form
\begin{equation}\label{Multivp}
	x\, '_{\boldsymbol g}(t)=  f(t,  x(t)),\quad x(t_0)= {x_0},
\end{equation}
with $t_0, T\in\mathbb R$, $T>0$, $X\subset \mathbb R^n$, $x_0\in X$ and $f:[t_0,t_0+T)\times X\to\mathbb R^n$, $f=(f_1,f_2,\dots,f_n)$. To that end, we introduce the concept of solution that is fundamental for the aims of this section. In what follows we denote by $I_\sigma=[t_0,t_0+\sigma)$, $\sigma\in(0,T]$, and $I=[t_0,t_0+T)$ and by $\overline I_\sigma$ and $\overline I$ the corresponding closure sets with respect to the usual topology in $\mathbb R$.
\begin{definition}\label{multivpsol}
	A \emph{solution}\index{solution of a $\boldsymbol g$--initial value problem} of \eqref{Multivp} on an interval $I_\sigma$, $\sigma\in (0,T]$, is a function $x\in\mathcal{AC}_{\boldsymbol g}(\overline I_\sigma, \mathbb R^n)$, $x=(x_1,x_2,\dots,x_n)$, such that $x(t_0)=x_0$, $x(t)\in X$ for all $t\in I_\sigma$ and
	\begin{equation}\label{eqmultivpsol}
		x'_{g_i}(t)=f_i(t,x(t)),\quad g_i\mbox{--a.a.} t\in I_\sigma,\quad i\in\{1,2,\dots,n\}.
	\end{equation}
	If $\sigma=T$, we say that $x$ is a \emph{global solution} of \eqref{Multivp}; otherwise, i.e. if $\sigma\in (0,T)$, we say that $x$ is a \emph{local solution} of \eqref{Multivp}.
\end{definition}

Before studying some existing and uniqueness results for \eqref{Multivp} we will reflect on the concept of classical solution in \cite{LoMa19}. When we talk about solutions in the classical sense, we mean solutions of \eqref{Multivp} that solve the problem on their whole interval of definition and have continuous derivatives. Of course, given the nature of Definition \ref{Stieltjesderivative}, this is impossible unless $C_{g_i}=\emptyset$, $i=1,2,\dots,n$. Nevertheless, we can talk about ``everywhere'' solutions referring to solutions solving the problem on the biggest set possible, i.e., excluding the corresponding set $C_{g_i}$. With this idea in mind, we start exploring some basic properties of ``everywhere'' solutions that will culminate in Theorem~\ref{Peanonosol}. The first results concern the continuity of the derivative.
\begin{proposition}\label{multigcompcor}
	Let $\sigma\in(0,T]$, $f:I\times\mathbb R^n\to\mathbb R^n$, $f=(f_1, f_2,\dots,f_n)$, be a $(\vec{\boldsymbol{g}}\times \Id)$--continuous function on $I_\sigma\times\mathbb R^n$ and $x\in\mathcal{AC}_{\boldsymbol g}(\overline I_\sigma,\mathbb R^n)$, $x=(x_1,x_2,\dots,x_n)$, be such that  $x(t_0)=x_0$ and
	\begin{equation}\label{solminusCgi}
	(x_i)'_{g_i}(t)=f_i(t,x(t)),\quad \mbox{for all }t\in I_\sigma\setminus C_{g_i},\quad i=1,2,\dots,n.
	\end{equation}
	Then, $(x_i)'_{g_i}$ is $\vec{\boldsymbol{g}}$--continuous on $I_\sigma\setminus C_{g_i}$.
\end{proposition}
\noindent
{\bf Proof.}
	Given that $x$ is $\boldsymbol g$--absolutely continuous on $\overline I_\sigma$,  we have that $x$ is $\boldsymbol g$--continuous on $\overline I_\sigma$, which implies that it is $\vec{\boldsymbol{g}}$--continuous on $\overline I_\sigma$. Thus Lemma \ref{pmbgcomp} yields that $f(\cdot,x(\cdot))$ is $\vec{\boldsymbol{g}}$--continuous. Hence, it follows that, for each $i\in\{1,2,\dots,n\}$, the map $f_i(\cdot,x(\cdot))$ is $\vec{\boldsymbol{g}}$--continuous and, since $(x_i)'_{g_i}$ is defined on $I_\sigma\setminus C_{g_i}$, the result follows.
\qed

It follows from Proposition \ref{multiconttocont} that, under the hypotheses of Proposition~\ref{multigcompcor}, $(x_i)'_{g_i}$ is $\widehat g$--continuous on $I_\sigma\setminus C_{g_i}$.
Furthermore, notice that if we replace the concept of $\vec{\boldsymbol{g}}$--continuity for $\boldsymbol g$--continuity in the hypotheses of Proposition \ref{multigcompcor} we still obtain $\vec{\boldsymbol{g}}$--continuous solutions. Nevertheless, with a similar reasoning to the one used for Proposition \ref{multigcompcor}, we can obtain the following result ensuring $\boldsymbol g$--continuity.

\begin{proposition}\label{multigcompcor2}
	Let $f:I\times\mathbb R^n\to\mathbb R^n$, $f=(f_1, f_2,\dots,f_n)$ and $x\in\mathcal{AC}_{\boldsymbol g}(\overline I_\sigma,\mathbb R^n)$, $x=(x_1,x_2,\dots,x_n)$, $\sigma\in (0,T]$, be such that  $x(t_0)=x_0$ and \eqref{solminusCgi} is satisfied.
	If there exists $i\in\{1,2,\dots,n\}$ such that $f_i(\cdot,x(\cdot))$ is $g_i$--continuous on $I_\sigma$, then $(x_i)'_{g_i}$ is $g_i$--continuous on $I_\sigma\setminus C_{g_i}$. In particular, if $f(\cdot,x(\cdot))$ is $\boldsymbol g$--continuous on $I_\sigma$, then $(x_i)'_{g_i}$ is $g_i$--continuous on $I_\sigma\setminus C_{g_i}$ for all $i\in\{1,2,\dots,n\}$.
\end{proposition}

Propositions \ref{multigcompcor} and \ref{multigcompcor2} show that, under suitable conditions, the solutions of \eqref{Multivp} have continuous derivatives in some sense. In particular, it is required that condition \eqref{solminusCgi} is satisfied. Nevertheless, solutions in the sense of Definition \ref{multivpsol} need not satisfy such condition. In the following results, we provide some conditions ensuring that \eqref{solminusCgi} or similar conditions are satisfied.

\begin{proposition}\label{multisoltoevery}
	Let $x=(x_1,x_2,\dots,x_n)$ be a solution of~\eqref{Multivp} on $I_\sigma$, $\sigma\in(0,T]$.
	Then:
	\begin{enumerate}
		\item[\textup{(i)}] If $f(\cdot,x(\cdot))$ is $\vec{\boldsymbol{g}}$--continuous on $I_\sigma$, then
		\[(x_i)'_{g_i}(t)=f_i(t,x(t)) \quad \mbox{for all $t \in (I_\sigma  \setminus (D_{\widehat g}\cup C_{g_i}))\cup D_{g_i}$,}\quad i=1,2,\dots,n.\]
		\item[\textup{(ii)}] If $f_i(\cdot,x(\cdot))$ is $g_i$--continuous for some $i\in\{1,2,\dots,n\}$, then
			\[(x_i)'_{g_i}(t)=f_i(t,x(t)) \quad \mbox{for all $t \in I_\sigma  \setminus C_{g_i}$}.\]
			In particular, if $f(\cdot,x(\cdot))$ is $\boldsymbol g$--continuous on $I_\sigma$, then \eqref{solminusCgi} holds.
	\end{enumerate}
\end{proposition}
\noindent
{\bf Proof.}
	Fix $i\in\{1,2,\dots,n\}.$ Since $  x$ is a solution of \eqref{Multivp}, we have that $  x_i\in\mathcal{AC}_{g_i}(\overline I_\sigma,\mathbb R^n)$ and $(x_i)'_{g_i}(t)=f_i(t,  x(t))$ for $g_i$--a.a. $t\in I_\sigma$.
	Hence, in particular, $(x_i)'_{g_i}(t)=f_i(t,  x(t))$ for all $t\in D_{g_i}$, so it is enough to show that the correponding equalities holds for all $t\in I_\sigma\setminus(C_{g_i}\cup D_{\widehat g})$ for (i), and for all $t\in I_\sigma\setminus(C_{g_i}\cup D_{g_i})$ for (ii). We will first prove (ii) and then, by making small modifications to that proof, we will obtain (i).

	Fix $t\in I_\sigma\setminus(C_{g_i}\cup D_{g_i}).$ Since $g_i$ is not constant on any neighbourhood of $t$, we may have $g_i(s)<g_i(t)$ for all $s<t$, $g_i(s)>g_i(t)$ for all $s>t$, or both. If $g_i(s)<g_i(t)$ for all $s<t$ and $t>t_0$, then, for all $s \in [t_0,t)$, the Fundamental Theorem of Calculus yields
	\begin{align*}
		\mu_{g_i}([s,t))\inf_{s\le r <t}f_i(r,  x(r))\le \int_{[s,t)}f_i(r,  x(r))  \dif g_i(r)=x_i(t)-x_i(s)\le\mu_{g_i}([s,t))\sup_{s\le r <t}f_i(r,  x(r)),
	\end{align*}
	which implies that
	\begin{equation}\label{gideraux}
		\inf_{s\le r<t}f_i(r,  x(r))\le\dfrac{x_i(s)-x_i(t)}{g_i(s)-g_i(t)}\le \sup_{s\le r<t}f_i(r,  x(r)).
	\end{equation}
	On the other hand, if $g_i(s)>g_i(t)$ for all $s>t$, then, following an analogous reasoning, we deduce that
	\begin{equation}\label{gideraux3}
		\inf_{t\le r <s}f_i(r,  x(r))\le\dfrac{x_i(s)-x_i(t)}{g_i(s)-g_i(t)}\le \sup_{t\le r <s}f_i(r,  x(r)).
	\end{equation}

	Now, for (ii), assume that $f_i(\cdot,x(\cdot))$ is $g_i$--continuous on $\overline I_\sigma$. In that case, since $g_i$ is continuous at $t$, we have that $f_i(\cdot,x(\cdot))$ is continuous at $t$. Therefore, if $t$ is such that $g_i(s)<g_i(t)$ for all $s<t$, \eqref{gideraux} implies that the following limit exists and
	\begin{equation}\label{gideraux2}
		\lim_{s \to t^-}\frac{x_i(s)-x_i(t)}{g_i(s)-g_i(t)}=f_i(t,  x(t)).
	\end{equation}
	If $g_i(s)=g_i(t)$ on some $[t,t+\delta]$, $\delta >0$, then the limit in \eqref{gideraux2} is $(x_i)'_{g_i}(t)$ and the proof is complete. Similarly, if $t$ is such that $g_i(s)>g_i(t)$ for all $s>t$, the continuity of $f_i(\cdot,x(\cdot))$ at $t$ and \eqref{gideraux3} ensure that the following limit exists and
	\begin{equation}\label{gideraux4}
		\lim_{s \to t^+}\frac{x_i(s)-x_i(t)}{g_i(s)-g_i(t)}=f_i(t,  x(t)).
	\end{equation}
	This covers all of the remaining cases, so the proof is finished for $f_i(\cdot,x(\cdot))$ $g_i$--continuous, and subsequently, for $f(\cdot,x(\cdot))$ $\boldsymbol g$--continuous.

	On the other hand, if $t\in I_\sigma\setminus(C_{g_i}\cup D_{\widehat g})$, we have that $t\in I_\sigma\setminus(C_{g_i}\cup D_{g_i})$. Hence, \eqref{gideraux} holds if $t$ is such that $g_i(s)<g_i(t)$ for all $s<t$, and \eqref{gideraux3}, if $g_i(s)>g_i(t)$ for all $s>t$. Now, under the hypotheses of (i), we have that $f(\cdot,x(\cdot))$ is $\vec{\boldsymbol{g}}$--continuous, so it is $\widehat g$--continuous at $t$. As a consequence, $f(\cdot,x(\cdot))$ is continuous at $t$ as $\widehat g$ is continuous at that point. Hence, by making analogous reasonings, we can obtain \eqref{gideraux2} and \eqref{gideraux4} to finish the proof.
\qed

Essentially, the proof of Proposition \ref{multisoltoevery} is a revision of the proofs of \cite[Proposition 7.6]{FriLo17} and \cite[Proposition 4.6]{LoMa19} in the context of \eqref{Multivp}. Of course, in the latter, the authors worked on the same framework as in this paper. Nevertheless, \cite[Proposition 4.6]{LoMa19} is not correct due to Definitions~\ref{multicont} and \ref{multigIdcont} not being equivalent to Definitions~\ref{vecmulticont} and \ref{vecgIdcont}, respectively. The following result serves as a proper reformulation of \cite[Proposition 4.6]{LoMa19} based on Proposition \ref{multisoltoevery}.

\begin{theorem} \label{multieverywheresolutions}
	Let $\sigma\in(0,T]$, $f:I\times\mathbb R^n\to\mathbb R^n$ be a $(\vec{\boldsymbol{g}}\times \Id)$--continuous function on $I_\sigma\times\mathbb R^n$ and $x=(x_1,x_2,\dots,x_n)$ be a solution of~\eqref{Multivp} on $I_\sigma$. If there exists $i\in\{1,2,\dots,n\}$ such that
	\begin{equation}\label{DhatgDg_i}
		(x_i)'_{g_i}(t)=f_i(t,x(t)),\quad \mbox{for all $t \in D_{\widehat g}\setminus D_{g_i}$,}
	\end{equation}
	then
	\begin{equation}\label{DhatgDg_i2}
		(x_i)'_{g_i}(t)=f_i(t,x(t)) \quad \mbox{for all $t \in I_\sigma  \setminus C_{g_i}$,}
	\end{equation}
	and $(x_i)'_{g_i}$ is $\vec{\boldsymbol{g}}$--continuous on $I_\sigma\setminus C_{g_i}$. In particular, if \eqref{DhatgDg_i} holds for all $i\in\{1,2,\dots,n\}$, then \eqref{DhatgDg_i2} holds for all $i\in\{1,2,\dots,n\}$ and each  $(x_i)'_{g_i}$, $i\in\{1,2,\dots,n\}$, is $\vec{\boldsymbol{g}}$--continuous on $I_\sigma\setminus C_{g_i}$.
\end{theorem}
\noindent
{\bf Proof.}
By definition, $x$ is $\boldsymbol g$--absolutely continuous on $\overline I_\sigma$ which implies that it is $\vec{\boldsymbol{g}}$--continuous on that set. Hence, given that $f$ is $(\vec{\boldsymbol{g}}\times \Id)$--continuous on $\overline I_\sigma\times\mathbb R^n$, Proposition \ref{multisoltoevery} ensures that
\[
(x_i)'_{g_i}(t)=f_i(t,x(t)) \quad \mbox{for all $t \in I_\sigma  \setminus (D_{\widehat g}\cup C_{g_i})\cup D_{g_i}$,}\quad i=1,2,\dots,n.
\]
Hence, for each $i\in\{1,2,\dots,n\}$ such that \eqref{DhatgDg_i} holds, we have that $(x_i)'_{g_i}(t)=f_i(t,x(t))$, $t\in I_\sigma\setminus C_{g_i}$. Now, Proposition \ref{multigcompcor} ensures that $(x_i)'_{g_i}$ is $\vec{\boldsymbol{g}}$--continuous in $I_\sigma\setminus C_{g_i}$, which finishes the proof. \qed

Observe that \eqref{DhatgDg_i} becomes vacuous when $D_{g_j}=D_{g_k}$ for all $j,k\in\{1,2,\dots,n\}$ which, as stated by \textup{Proposition \ref{equivpmbgvecg}}, is guaranteed to happen when \textup{Definitions~\ref{multicont}} and \textup{\ref{vecmulticont}} are equivalent. This justifies the statement of \textup{\cite[Proposition 4.6]{LoMa19}}, since the authors wrongly used both definitions equivalently. This, of course, affected \textup{\cite[Theorem 4.8]{LoMa19}}, where the authors guaranteed the existence of a local solution (and everything that follows from \textup{\cite[Proposition 4.6]{LoMa19}}) under the assumption of an extra hypothesis: for each $i\in\{1,2,\dots,n\}$, there exists $h_i\in\mathcal L^1_{g_i}(I,[0,+\infty))$ such that
\begin{equation}\label{extracondPeano}
	|f_i(t,x)|\le h_i(t),\quad g_i\mbox{--a.a. $t\in I$},\quad x\in \overline{B(x_0, r)}.
\end{equation}
It is possible to obtain a correct formulation of \textup{\cite[Theorem 4.8]{LoMa19}} by noting that the $(\vec{\boldsymbol{g}}\times \Id)$--continuity of $f$ together with condition \eqref{extracondPeano} are enough to obtain the existence of a local solution through \textup{\cite[Theorem 4.5]{LoMa19}}. After that, all that is left to do is to consider \textup{Theorem \ref{multieverywheresolutions}} to obtain the right version of \textup{\cite[Theorem 4.8]{LoMa19}}.

As a final note, we obtain an analogous result to Theorem \ref{multieverywheresolutions} yielding $\boldsymbol g$--continuity instead of $\vec{\boldsymbol{g}}$--continuity. This result follows from Proposition \ref{multigcompcor2} combined with statement (ii) in Proposition~\ref{multisoltoevery}.

\begin{theorem}\label{Peanonosol}
Let $x=(x_1,x_2,\dots,x_n)$ be a solution of~\eqref{Multivp} on $I_\sigma$, $\sigma\in (0,T]$. If $f_i(\cdot, x(\cdot))$ is $g_i$--continuous on $I_\sigma$ for some $i\in\{1,2,\dots,n\}$, then
\begin{equation*}
		(x_i)'_{g_i}(t)=f_i(t,x(t)) \quad \mbox{for all $t \in I_\sigma  \setminus C_{g_i}$,}
\end{equation*}
and $(x_i)'_{g_i}$ is $g_i$--continuous on $I_\sigma\setminus C_{g_i}$. In particular, if $f(\cdot, x(\cdot))$ is $\boldsymbol g$--continuous on $I_\sigma$, then \eqref{solminusCgi} holds and $(x_i)'_{g_i}$ is $g_i$--continuous on $I_\sigma\setminus C_{g_i}$ for all $i\in\{1,2,\dots,n\}$.
\end{theorem}

\subsection{Uniqueness of solution}
We now continue the work in \cite{LoMa19} by researching some uniqueness conditions for \eqref{Multivp}. In particular, we shall adapt the results in \cite{MaMon20} to the context of differential equations with several Stieltjes derivatives. For this endeavour, as well as for the question of existence of solution, we can assume that $\boldsymbol g$ is continuous at $t_0$, as pointed out in \cite{LoMa19}. To see that this is the case, it is enough to use an analogous reasoning to that in \cite[Section 5]{FriLo17}.

Our first uniqueness criterion for \eqref{Multivp} is the analogous to \cite[Theorem 3.9]{MaMon20} and it is inspired by the ideas of \cite[Theorem 4.8]{Schw92}. In order to obtain the mentioned result we need the following reformulation of \cite[Theorem 1.40]{Schw92} in the context of the Lebesgue--Stieltjes integral, see
\cite[Lemma~3.8]{MaMon20}.
\begin{lemma}\label{lem}
	Let $h:\mathbb R\to\R$ be a nondecreasing left--continuous function, and let $\omega:[0,+\infty)\to[0,+\infty)$ be a continuous nondecreasing function such that $\omega(0)=0$, $\omega(s)>0$ for $s>0$. For a fixed $u_0>0$, define
	\[\Omega(r)=\int_{u_0}^r\frac{1}{\omega(s)}\dif s,\quad r\in(0,+\infty).\]
	and denote by $\alpha=\lim_{r\to0^+}\Omega(r)\ge-\infty$, $\beta=\lim_{r\to+\infty}\Omega(r)\le+\infty$.
	If $\psi:[a,b]\to[0,+\infty)$ is a bounded function, and there exists $\kappa>0$ such that
	\begin{equation*}\label{intineq}
	\psi(s)\le \kappa+\int_{[a,s)}\omega(\psi(\sigma))\dif h(\sigma),\quad s\in [a,b],
	\end{equation*}
	and $\Omega(\kappa)+h(b)-h(a)<\beta$, then
	\[\psi(s)\le\Omega^{-1}(\Omega(\kappa)+h(s)-h(a)),\quad s\in [a,b],\]
	where $\Omega^{-1}:(\alpha,\beta)\to\mathbb R$ is the inverse function of $\Omega$.
\end{lemma}

We are now able to state and prove the following uniqueness result under an Osgood type condition.

\begin{theorem}\label{multiosgood}
	Let $X\subset\mathbb R^n$ be such that $x_0\in X$, $f:[t_0,t_0+T)\times X\to\mathbb R^n$, $f=(f_1,f_2,\dots,f_n)$, and $\omega:[0,+\infty)\to[0,+\infty)$ be a nondecreasing continuous function such that $\omega(0)=0$, $\omega(s)>0$ for all $s>0$ and
	\begin{equation}\label{osgoodcond}
	\lim_{\varepsilon\to0^+}\int_{\varepsilon}^{u_0}\frac{1}{\omega(s)}\dif s=+\infty,
	\end{equation}
	for some $u_0>0$.
	If there exists $\sigma\in (0,T]$ such that for each $i\in\{1,2,\dots,n\}$,
	\[|f_i(t,x)-f_i(t,y)|\le \omega(\|x-y\|),\quad g_i\mbox{--a.a. }t\in I_\sigma,\quad x,y\in X,\]
	then \eqref{Multivp} has at most one solution on $I_\sigma$.
\end{theorem}
\noindent
{\bf Proof.}
		Let $x$, $y$ be solutions of \eqref{Multivp} on $I_\sigma$.
		Define $\psi:\overline I_\sigma\to[0,+\infty)$, $\Omega:(0,+\infty)\to(0,+\infty)$ as
		\begin{equation}\label{psiOmega}
		\psi(t)=\|x(t)-y(t)\|,\quad\quad\quad\Omega(r)=\displaystyle\int_{u_0}^r \frac{1}{\omega(s)}\dif s.
		\end{equation}

		First, note that,
		given that $x,y\in\mathcal{AC}_{\boldsymbol g}(\overline I_\sigma,\mathbb R^n)$, we have that each component of $x-y$ is Borel measurable (see Remark \ref{multigmeasurable}) and, since $\psi$ is the pointwise maximum of Borel measurable maps, we have that $\psi$ is Borel measurable. Now, given that $\omega$ is continuous, it follows that $\omega\circ \psi$ is Borel measurable, which guarantees that it is $\widehat g$ and $g_i$--measurable, $i=1,2,\dots,n$. Moreover, Remark~\ref{gabsbound} ensures that $x-y$ is bounded, which implies that so is $\psi$, yielding that $\omega\circ \psi$ is bounded as well. Hence, it follows that $\omega\circ \psi$ is integrable with respect to $\widehat g$ and $g_i$, $i=1,2,\dots,n$.

		Let $K>0$ be an upper bound of $\omega\circ \psi$. Without loss of generality, we assume that $\boldsymbol g$ is continuous, which ensures that $\widehat g$ is continuous at $t_0$. Then, for each $\gamma\in (0,\sigma)$,
		\[\int_{[t_0,t_0+\gamma)}\omega(\psi(s))\dif \widehat g(s)\le\int_{[t_0,t_0+\gamma)}K\dif \widehat g(s)=K\mu_{\widehat g}([t_0,t_0+\gamma))<\widehat \varepsilon(\gamma),\]
		where $\widehat \varepsilon(\gamma)=K\mu_{\widehat g}([t_0,t_0+\gamma))+\gamma>0$.
		Noting that $\widehat \varepsilon$ and $\omega$ are in the same circunstances as $\varepsilon$ and $\omega$ in the proof of \cite[Theorem 3.9]{MaMon20} with $\widehat g$ in place of $g$,  we can repeat the same arguments there to see that there exists $0<R<\gamma$ such that
		\[\Omega(\widehat \varepsilon(\delta))+\widehat g(t_0+\sigma)-\widehat g(t_0+\delta)<\beta:=\lim_{r\to\infty}\Omega(r)\quad\mbox{for \ }\delta\in(0,R).\]
		By definition we have that, for each $t\in\overline I_\sigma$, there is $j_t\in\{1,2,\dots,n\}$ such that $\psi(t)=|x_{j_t}(t)-y_{j_t}(t)|$. Therefore, Theorem \ref{gFTC2} yields that for each $t\in\overline I_\sigma$,
		\begin{align*}
		\psi (t)&\le \int_{[t_0,t)}\left|f_{j_t}(s,x(s))-f_{j_t}(s,y(s))\right|\dif g_{j_t}(s)\le \int_{[t_0,t)}\omega(\|x(s)-y(s)\|)\dif g_{j_t}(s)\\
		&\le \sum_{i=1}^n\int_{[t_0,t)}\omega(\|x(s)-y(s)\|)\dif g_{i}(s)= \int_{[t_0,t)}\omega(\|x(s)-y(s)\|)\dif \widehat g(s)\\
		&=\int_{[t_0,t_0+\delta)}\omega(\psi(s))\dif \widehat g(s)+\int_{[t_0+\delta,t)}\omega(\psi(s))\dif \widehat g(s)<\widehat \varepsilon(\delta)+\int_{[t_0+\delta,t)}\omega(\psi(s))\dif \widehat g(s)&&
		\end{align*}
		for all $\delta\in(0,R)$. Therefore, the assumptions of Lemma~\ref{lem} are satisfied, which guarantees that
		\[\psi(t)\le \Omega^{-1}(\Omega(\varepsilon(\delta))+\widehat g(t)-\widehat g(t_0+\delta)),\quad\delta\in(0,R),\quad t\in \overline I_\sigma.\]
		Applying $\Omega$ on both sides of the inequality, we obtain
		\[\Omega(\psi(t))-\Omega(\varepsilon(\delta))\le \widehat g(t)-\widehat g(t_0+\delta)\le \widehat g(t)-\widehat g(t_0),\quad\delta\in(0,R),\quad t\in \overline I_\sigma.\]

		Suppose $\psi\not= 0$ on $\overline I_\sigma$. In that case, there must exist $t^*\in \overline I_\sigma$ such that $\psi(t^*)>0$. Then, for all $\delta\in(0,R)$ such that $\delta<t^*-t_0$,
		\[\int_{\varepsilon(\delta)}^{\psi(t^*)}\frac{1}{\omega(s)}\dif s=\Omega(\psi(t^*))-\Omega(\varepsilon(\delta))<\widehat g(t^*)-\widehat g(t_0),\]
		and, by taking the limit as $\delta\to 0^+$, we obtain
		\[\lim_{\delta\to0^+}\int_{\varepsilon(\delta)}^{\psi(t^*)}\frac{1}{\omega(s)}\dif s<\widehat g(t^*)-\widehat g(t_0)<+\infty,\]
		which contradicts (ii). Hence, we must have  $\psi=0$ on $\overline I_\sigma$, i.e. $x=y$ on that interval.\qed

		Observe that Theorem \ref{multiosgood} returns \cite[Theorem 3.9]{MaMon20} in the corresponding setting (namely, when $\boldsymbol g=(g,g,\dots,g)$ for some derivator $g$) as, for any $u\in(0,+\infty)$,
		\begin{equation}\label{osgoodcondsimpl}
		\lim_{\varepsilon\to0^+}\int_{\varepsilon}^{u}\frac{1}{\omega(s)}\dif s=\lim_{\varepsilon\to0^+}\int_{\varepsilon}^{u_0}\frac{1}{\omega(s)}\dif s+\int_{u_0}^{u}\frac{1}{\omega(s)}\dif s.
		\end{equation}
		Furthermore, it is possible to obtain the following more general result, which is the analog of the Montel--Tonelli uniqueness result in \cite{MaMon20} in the setting of \eqref{Multivp}. In its proof, we make use of the Kurzweil--Stieltjes integral, which we denote by
		\[\KS\int_{a}^b F(s) \dif g(s).\]
		For more information on this integral, we refer the reader to \cite{MonSlaTvr18}.

\begin{theorem}\label{multiMontel}
	Let $X\subset\mathbb R^n$ be such that $x_0\in X$, $f:[t_0,t_0+T)\times X\to\mathbb R^n$, $f=(f_1,f_2,\dots,f_n)$, and $\omega:[0,+\infty)\to[0,+\infty)$ be a nondecreasing continuous function such that $\omega(0)=0$, $\omega(s)>0$ for all $s>0$ and \eqref{osgoodcond} holds for some $u_0>0$.
	If there exists $\sigma\in (0,T]$ and $\varphi:I_\sigma\to[0,+\infty)$ such that for each $i\in\{1,2,\dots, n\}$, $\varphi\in \mathcal L^1_{g_i}(I_\sigma,[0,+\infty))$ and
	\begin{equation}\label{MontelTonellicond}
	|f_i(t,x)-f_i(t,y)|\le \varphi(t)\omega(\|x-y\|),\quad g_i\mbox{--a.a. }t\in I_\sigma,\quad x,y\in X,
	\end{equation}
	then \eqref{Multivp} has at most one solution on $I_\sigma$.
\end{theorem}
\noindent
{\bf Proof.}
		Let $x$, $y$ be solutions of \eqref{Multivp} on $I_\sigma$ and define $\psi:\overline I_\sigma\to[0,+\infty)$ as in \eqref{psiOmega}.
		We can show, using the same arguments as in Theorem \ref{multiosgood}, that $\omega\circ \psi$ is bounded and Borel measurable. Therefore, given $K>0$ an upper bound of $\omega\circ \psi$, we have that $|\omega(\psi(t))\varphi(t)|\le K \varphi(t) $, $t\in \overline I_\sigma$.
		This fact, together with Proposition \ref{gsum}, ensure that $\varphi\cdot\omega\circ\psi$ is $\widehat g$--integrable on $I_\sigma$.
		Define $\overline g:\mathbb R\to\mathbb R$ as
		\[
		\overline g(t)=\left\{
		\begin{array}{ll}
		0,& \mbox{\textup{if} }t\le t_0,\vspace{0.1cm}\\
		\displaystyle \int_{[t_0,t)}\varphi(s)\dif \widehat g(s),&\mbox{\textup{if} }\displaystyle t_0<t\le t_0+\sigma,\vspace{0.1cm}\\
		\displaystyle \int_{[t_0,t_0+\sigma)}\varphi(s)\dif \widehat g(s),&\mbox{\textup{if} }\displaystyle t>t_0+\sigma.
		\end{array}
		\right.
		\]
		Recalling the relation between the Lebesgue--Stieltjes and Kurzweil--Stieltjes integrals (see \cite[Chapter~6, Section~12]{MonSlaTvr18}) and the substitution formula for Kurzweil--Stieltjes integral, for each $t\in \overline I_\sigma$ we have
		\begin{equation}\label{multisubs}\begin{aligned}
		\int_{[t_0,t)}\omega(\psi(s))\varphi(s) \dif \widehat g(s) = &\KS\int_{t_0}^t\omega(\psi(s))\varphi(s) \dif \widehat g(s)\\ =& \KS\int_{t_0}^t\omega(\psi(s)) \dif \overline{g}(s)=\int_{[t_0,t)}\omega(\psi(s)) \dif \overline{g}(s).
		\end{aligned}
		\end{equation}
		Observe that the integrals in \eqref{multisubs} are well--defined in the Lebesgue--Stieltjes sense as the map $\omega\circ \psi$ is Borel measurable. Once again, without loss of generality, we assume that $\boldsymbol g$ is continuous at $t_0$, which implies that so is $\widehat{g}$ and, as a consequence, so is $\overline g$.
		This, together with \eqref{multisubs}, yields
		\begin{equation*}\label{e-2}
		\int_{[t_0,t_0+\gamma)}\omega(\psi(s)) \dif \overline{g}(s)
		\le\int_{[t_0,t_0+\gamma)}K\dif \overline{g}(s)<K\mu_{\overline{g}}([t_0,t_0+\gamma))+\gamma=:\overline\varepsilon(\gamma),\quad \gamma\in(0,\sigma).
		\end{equation*}
		Therefore, the result can be proved by reasoning as in the proof of Theorem \ref{multiosgood} with the appropriate adjustments, i.e. replacing $\widehat g$ by $\overline{g}$ and $\widehat\varepsilon(\cdot)$ by $\overline\varepsilon(\cdot)$ accordingly.\qed

\subsection{Existence and uniqueness of solution}
As a final step on the study of \eqref{Multivp}, we combine the uniqueness results in the previous section with some information available in the literature regarding the existence of solution. Specifically, we will use \cite[Theorem 4.5]{LoMa19}, which requires the following definition, a slightly more general version of \cite[Definition 4.7]{FriLo17}.

\begin{definition}\label{gfCarat}
	Let  $g:\mathbb R\to\mathbb R$ be a nondecreasing and left--continuous function, $J\subset \mathbb R$ and $X\subset \mathbb R^m$, $X\not=\emptyset$.
	A function $f:J\times X\to\mathbb R^n$ is said to be \emph{$g$--Carath\'eodory} if the following properties are satisfied:
	\begin{enumerate}
		\item[\textup{(i)}] For each $x\in X$, the map $f(\cdot,x)$ is $g$--measurable.
		\item[\textup{(ii)}] The map $f(t,\cdot)$ is continuous for $g$--a.a. $t\in J$.
		\item[\textup{(iii)}] For each $r>0$, there exists $h_r\in \mathcal L^1_g(J,[0,+\infty))$ such that
		\[\|f(t,x)\|\le h_r(t),\quad \mbox{$g$--a.a. $t\in J$,}\quad \mbox{$x\in X$, $\|x\|\le r$.}\]
	\end{enumerate}
\end{definition}
\begin{remark}\label{remosgoodcarat}
	As proved in \textup{\cite[Theorem 4.1]{MaMon20}}, a sufficient condition for the $f$ in \eqref{Multivp} to be $g$--Carath\'eodory for some derivator $g$ is that $f(\cdot, x_0)\in\mathcal L^1_g(I_\sigma,\mathbb R)$ and conditions \textup{(i)} in \textup{Definition~\ref{gfCarat}} and \eqref{MontelTonellicond} with $\omega$ continuous at $0$, $\omega(0)=0$, are satisfied.
\end{remark}

Based on this definition, we find the following Peano--type existence result in \cite{LoMa19}.
\begin{theorem} \label{multipeano}
	Let $r>0$ and $f:I\times\overline{B( {x_0},r)}\to\mathbb R^n$ be such that for every $i=1,2,\dots,n$, the map $f_i$ is $g_i$--Carath\'eodory. Then, there exists $\sigma\in(0,T]$ such that \eqref{Multivp} has a solution on $I_\sigma$.
\end{theorem}

As a direct consequence of Theorem \ref{multipeano} and Remark \ref{remosgoodcarat}, we obtain the following Montel--Tonelli existence result in a similar fashion to \cite[Theorem 4.2]{MaMon20}. In fact, Theorem \ref{multiexistence} reduces to the mentioned result in the corresponding setting, i.e., for $\boldsymbol g=(g,g,\dots,g)$ for a derivator $g$.
\begin{theorem}\label{multiexistence}
	Let $r>0$ and $f:I\times \overline{B(x_0,r)}\to\mathbb R^n$ be such that, for each $i\in\{1,2,\dots,n\}$, the following conditions are satisfied:
	\begin{enumerate}
		\item[\textup{(i)}] For each $x\in \overline{B(x_0,r)}$, $f_i(\cdot, x)$ is $g_i$--measurable.
		\item[\textup{(ii)}] $f_i(\cdot, x_0)\in\mathcal L^1_{g_i}(I,\mathbb R)$.
		\item[\textup{(iii)}] There exist $\varphi_i\in \mathcal L^1_{g_i}(I,[0,+\infty))$ and $\omega_i:[0,+\infty)\to[0,+\infty)$ nondecreasing, continuous at $0$ with $\omega_i(0)=0$ and such that \eqref{MontelTonellicond} holds for $\sigma=T$ and $X=\overline{B(x_0,r)}$.
	\end{enumerate}
	Then, there exists $\sigma\in(0,T]$ such that \eqref{Multivp} has a solution on $I_\sigma$.
\end{theorem}

Naturally, combining the hypotheses of Theorems \ref{multiMontel} and \ref{multiexistence} we can obtain an Montel--Osgood--Tonelli type existence and uniqueness result of local solutions of \eqref{Multivp}. This result is a generalization of \cite[Theorem 4.3]{MaMon20}. To see that this is the case, it is enough to bear in mind expression \eqref{osgoodcondsimpl}.

\begin{theorem}\label{multiexist-uniq}
	Let $r>0$,
	$f:I\times \overline{B( {x_0},r)}\to\mathbb R^n$ and $\omega:[0,+\infty)\to[0,+\infty)$ be a nondecreasing continuous function such that $\omega(0)=0$ and $\omega(s)>0$, $s>0$. Assume that the following conditions are satisfied:
	\begin{enumerate}
		\item[\textup{(i)}] For every $i=1,2,\dots, n$ and each $x\in \overline{B(x_0,r)}$, the map $f_i(\cdot, x)$ is $g_i$--measurable.
		\item[\textup{(ii)}] For every $i=1,2,\dots, n$, $f_i(\cdot, x_0)\in\mathcal L^1_{g_i}(I,\mathbb R)$.
		\item[\textup{(iii)}] There exists $u_0>0$ such that \eqref{osgoodcond} holds.
		\item[\textup{(iv)}] There exists a map
		$\varphi:I\to[0,+\infty)$ such that, for each $i\in\{1,2,\dots, n\}$, $\varphi\in \mathcal L^1_{g_i}(I,[0,+\infty))$ and
		\[|f_i(t,x)-f_i(t,y)|\le \varphi(t)\omega(\|x-y\|),\quad g_i\mbox{--a.a. }t\in I,\quad x,y\in \overline{B(x_0,r)}.\]
	\end{enumerate}
	Then, there exists $\sigma\in(0,T]$ such that \eqref{Multivp} has a unique solution on $I_\sigma$.
\end{theorem}

Observe that for the particular choice of $\omega(r)=r$, $r\ge 0$, we obtain \cite[Theorem 4.4]{LoMa19}. As commented in Section~\ref{sectionintegral}, the proof of such result relies, implicitly, on the fact that if $\varphi\in \mathcal L^1_{g_i}(I,[0,+\infty))$ for all $i\in\{1,2,\dots,n\}$, then $\varphi\in \mathcal L^1_{\widehat g}(I,[0,+\infty))$, which was unclear in that setting. Here, and specifically, in the proof of Theorem \ref{multiMontel} we make use of Proposition \ref{gsum} to ensure that this is the case, hence improving the work in \cite{LoMa19}.

Note that Theorem \ref{multiexist-uniq} is stated in the context of a neighborhood of the $x_0$. Following the steps of \cite{MaMon20}, we can obtain an Montel--Osgood--Tonelli type existence and uniqueness result with the map $f$ in \eqref{Multivp} defined on the whole $I\times \mathbb R^n$. To that end, we introduce the following lemma, which is a reformulation of \cite[Theorem 4.4]{MaMon20} in the context of \eqref{Multivp}.

\begin{lemma}\label{multibound-sol}
	Let $f : I \times \mathbb R^n \to {\mathbb R}^n$ and $\omega:[0,+\infty)\to[0,+\infty)$ be a nondecreasing continuous function such that $\omega(0)=0$ and $\omega(s)>0$, $s>0$.
	Assume that $\boldsymbol g$ is continuous at $t_0$ and that the following conditions are satisfied:
	\begin{enumerate}
		\item[\textup{(i)}] For every $i=1,2,\dots, n$ and $  x\in \mathbb R^n$, the map $f_i(\cdot,  x)$ is $g_i$--measurable.
		\item[\textup{(ii)}] For each $i\in\{1,2,\dots,n\}$, $f_i(\cdot, x_0)\in\mathcal L^1_{g_i}(I,\mathbb R)$.
		\item[\textup{(iii)}] There exists $u_0>0$ such that \eqref{osgoodcond} holds.
		\item[\textup{(iv)}] There exists a map
		$\varphi:I\to[0,+\infty)$ such that, for each $i\in\{1,2,\dots, n\}$, $\varphi\in \mathcal L^1_{g_i}(I,[0,+\infty))$ and
		\[|f_i(t,x)-f_i(t,y)\|\le \varphi(t)\omega(\|x-y\|),\quad g_i\mbox{--a.a. }t\in I,\quad x,y\in  \mathbb R^n.\]
	\end{enumerate}
	Then there exist $t_1\in(t_0,t_0+T]$ and a nondecreasing function $h:[t_0,t_1]\to\R$ such that for every solution of \eqref{Multivp}, $x:\overline I_\sigma\to\R^n$, $\sigma\in(0,T]$, we have
	\[
	\|x(t)-x_0\|\le h(t),\quad t\in \overline I_\sigma\cap [t_0,t_1].
	\]
\end{lemma}
\noindent
{\bf Proof.}
	Define $\kappa:\overline I\to\mathbb R$ as
	\[
	\kappa(t)=\int_{[t_0,t)}|f_i(s,x_0)|\dif \widehat{g}(s),\quad t\in \overline I.
	\]
	Note that hypothesis (ii) and Proposition \ref{gsum} ensure that $\kappa$ is well--defined.
	Let $x:\overline I_\sigma\to\R^n$ be a solution of \eqref{Multivp}. Remark~\ref{pmbgabstohatgabs} ensures that $x-x_0\in\mathcal{AC}_{\widehat g}(\overline I_\sigma,\mathbb R^n)$, which yields $\|x-x_0\|\in\mathcal{AC}_{\widehat g}(\overline I_\sigma,\mathbb R^n)$. Hence, it is Borel measurable and bounded. Given that $\omega$ is continuous, it follows that $\omega(\|x-x_0\|)$ is Borel measurable and thus, $\widehat g$ and $g_i$--measurable for each $i\in\{1,2,\dots,n\}$. Furthermore, since $\omega$ is nondecreasing, it follows that it is bounded, which guarantees that the map $\varphi(t)\omega(\|x(t)-x_0\|)$ is $\widehat g$ and $g_i$--integrable for each $i\in\{1,2,\dots,n\}$.

	Now, for each $t\in\overline I_\sigma$, there is $j_t\in\{1,2,\dots,n\}$ such that $\|x(t)-x_0\|=|x_{j_t}(t)-x_{0,j_t}(t)|$. Hence, condition (iv) yields that, for each $t\in\overline I_\sigma$,
	\begin{align*}
	\|x(t)-x_0\|&\le \int_{[t_0,t)}\left|f_{j_t}(s,x_0)\right|\dif g_{j_t}(s)+\int_{[t_0,t)}\left|f_{j_t}(s,x(s))-f_{j_t}(s,x_0)\right|\dif g_{j_t}(s)\\
	&\le  \sum_{i=1}^n \left(\int_{[t_0,t)}|f_i(s,x_0)|\dif g_i(s)+\int_{[t_0,t)}\varphi(s)\omega(\|x(s)-x_0\|)\dif g_{i}(s)\right)\\
	&\le \kappa(t)+\int_{[t_0,t)}\varphi(s)\omega(\|x(s)-x_0\|)\dif \widehat g(s).
	\end{align*}
	Define $\overline g:\mathbb R\to\mathbb R$ as
	\[
	\overline g(t)=\left\{
	\begin{array}{ll}
	0,& \mbox{\textup{if} }t\le t_0,\vspace{0.1cm}\\
	\displaystyle \int_{[t_0,t)}\varphi(s)\dif \widehat g(s),&\mbox{\textup{if} }\displaystyle t_0<t\le t_0+\sigma,\vspace{0.1cm}\\
	\displaystyle \int_{[t_0,t_0+\sigma)}\varphi(s)\dif \widehat g(s),&\mbox{\textup{if} }\displaystyle t>t_0+\sigma.
	\end{array}
	\right.
	\]
	Recalling the relation \eqref{multisubs}, it follows that
	\begin{equation}\label{multisol-1}
	\|x(t)-x_0\|\le \kappa(t)+\int_{[t_0,t)}\omega(\|x(s)-x_0\|)\dif \overline{g}(s),\quad t\in \overline I_\sigma.
	\end{equation}
	In order to apply \cite[Theorem 1.40]{Schw92}, consider the map $\Omega:(0,+\infty)\to(0,+\infty)$ as in \eqref{psiOmega}.
	Since $\lim_{r\to0^+}\Omega(r)=-\infty$, there exists $R>0$ such that
	\[
	\Omega(R)+\overline{g}(t_0+T)-\overline{g}(t_0)<\beta:=\lim_{r\to\infty}\Omega(r)\le +\infty.
	\]
	Since $\boldsymbol g$ is continuous at $t_0$, we have that  $\widehat g$ is continuous at $t_0$ as well. Then, we can find $t_1\in(t_0,t_0+T]$ such that $\kappa(t_1)\le R$. The monotonicity of $\Omega$ then yields
	\[
	\Omega(\kappa(t_1))+\overline{g}(t_1)-\overline{g}(t_0)<\beta.
	\]
	The inequality above together with \eqref{multisol-1} shows that the assumptions of Lemma~\ref{lem} are satisfied on the interval $[t_0,t_1]$, therefore
	\[
	\|x(t)-x_0\|\le \Omega^{-1}(\Omega(\kappa(t_1))+\overline{g}(t)-\overline{g}(t_0))=:h(t),\quad t\in \overline I_\sigma\cap [t_0,t_1],
	\]
	and $h:[t_0,t_1]\to\R$ is the desired monotone function.
\qed

Furthermore, for the proof of Theorem~\ref{multiexist-Schauder}, we will also require the following result which is a reformulation of \cite[Proposition~3.5]{LoMa19} based on \cite[Proposition~5.5]{FriLo17} and Remark~\ref{gabsbound}.
\begin{proposition}\label{multirelcomp}
	Let $\mathcal S$ be a subset of $\mathcal{AC}_{\boldsymbol g}([a,b],\mathbb R^n)$. Assume that for each $i\in\{1,2,\dots,n\}$ the following conditions are satisfied:
	\begin{enumerate}
		\item[\textup{(i)}] The set $\{F_i(a):   F=(F_1,F_2,\dots,F_n)\in \mathcal S\}$ is bounded.
		\item[\textup{(ii)}] There exists $h_i\in\mathcal L_{g_i}^1([a,b),[0,+\infty))$ such that
		\[|(F_i)'_{g_i}(t)|\le h_i(t),\quad g_i\mbox{--a.a. }t\in[a,b),\quad \mbox{for all }  F=(F_1,F_2,\dots,F_n)\in \mathcal S.\]
	\end{enumerate}
	Then $\mathcal S$ is a relatively compact subset of $\mathcal{BC}_{\boldsymbol g}([a,b],\mathbb R^n)$.
\end{proposition}

Now, we can state and prove the following Montel--Osgood--Tonelli
existence and uniqueness of solution result for problem \eqref{Multivp}. Observe that, although we are imposing global conditions on $\mathbb R^n$, we can only ensure the existence of a unique local solution. In that sense, Theorem \ref{multiexist-Schauder} is only a partial improvement with respect to \cite[Theorem 4.3]{LoMa19}.
\begin{theorem}\label{multiexist-Schauder}
	Let $f : I \times \mathbb R^n \to {\mathbb R}^n$ and $\omega:[0,+\infty)\to[0,+\infty)$ be a nondecreasing continuous function such that $\omega(0)=0$ and $\omega(s)>0$, $s>0$.
	Assume that the following conditions are satisfied:
	\begin{enumerate}
		\item[\textup{(i)}] For every $i=1,2,\dots, n$ and $  x\in \mathbb R^n$, the map $f_i(\cdot,  x)$ is $g_i$--measurable.
		\item[\textup{(ii)}] For each $i\in\{1,2,\dots,n\}$, $f_i(\cdot, x_0)\in\mathcal L^1_{g_i}(I,\mathbb R)$.
		\item[\textup{(iii)}] There exists $u_0>0$ such that \eqref{osgoodcond} holds.
		\item[\textup{(iv)}] There exists a map
		$\varphi:I\to[0,+\infty)$ such that, for each $i\in\{1,2,\dots, n\}$, $\varphi\in \mathcal L^1_{g_i}(I,[0,+\infty))$ and
		\[|f_i(t,x)-f_i(t,y)|\le \varphi(t)\omega(\|x-y\|),\quad g_i\mbox{--a.a. }t\in I,\quad x,y\in  \mathbb R^n.\]
	\end{enumerate}
	Then there exists $\sigma\in(0,T]$ such that \eqref{Multivp} has a unique solution on $I_\sigma$.
\end{theorem}
\noindent
{\bf Proof.}
	Without loss of generality, we assume that $\boldsymbol g$ is continuous at $t_0$. Let $h:[t_0,t_1]\to \R$ be the function whose existence is guaranteed by Lemma \ref{multibound-sol} and denote $R=h(t_1)$. Define $\overline g:\mathbb R\to\mathbb R$ as
	\[
	\overline g(t)=\left\{
	\begin{array}{ll}
	0,& \mbox{\textup{if} }t\le t_0,\vspace{0.1cm}\\
	\displaystyle \int_{[t_0,t)}\varphi(s)\dif \widehat g(s),&\mbox{\textup{if} }\displaystyle t_0<t\le t_0+T,\vspace{0.1cm}\\
	\displaystyle \int_{[t_0,t_0+T)}\varphi(s)\dif \widehat g(s),&\mbox{\textup{if} }\displaystyle t>t_0+T.
	\end{array}
	\right.
	\]
	Since $\boldsymbol g$ is continuous at $t_0$, we have that so are $g_i$, $i=1,2,\dots,n$, and $\widehat{g}$ and, as a consequence, $\overline g$ is continuous at $t_0$ as well. Thus, we can choose $\sigma\in(0,T]$ such that $t_0+\sigma\le t_1$ and
	\begin{equation}\label{B-R}
	\omega(R)\mu_{\overline{g}}([t_0,t_0+\sigma))+\sum_{i=1}^n \int_{[t_0,\sigma)}|f_i(s,x_0)|\dif g_i(s)<R.
	\end{equation}

	Consider $B=\big\{x\in \mathcal{BC}_{\boldsymbol g}(\overline I_\sigma,\R^n):\|x(t)-x_0\|\le R,\,t\in \overline I_\sigma\big\}$. Clearly, $B$ is a closed and convex subset of $\mathcal{BC}_{\boldsymbol g}(\overline I_\sigma,\R^n)$. Now let us define $F:B\to \mathcal{BC}_g(\overline I_\sigma,\R^n)$ by
	\[
	Fx(t)=x_0+\int_{[t_0,t)}f(s,x(s))\dif \boldsymbol g(s),\quad t\in I_\sigma.
	\]
	It follows from Remark~\ref{remosgoodcarat} that $f_i$ is $g_i$--Carath\'eodory, $i=1,2,\dots,n$. Furthermore, given $x\in B$, we have that $x$ is Borel measurable and therefore, $g_i$--measurable for each $i\in\{1,2,\dots,n\}$. Thus, \cite[Lemma~7.2]{FriLo17} and Remark \ref{gFTC} ensure that $F$ is well-defined. Moreover, the continuity of $\omega$ together with condition (iv) implies that $F$ is continuous. Furthermore, following a similar argument as the one used to obtain \eqref{multisol-1}, and using the relation \eqref{multisubs}, it follows from condition (iv) and \eqref{B-R} that for $x\in B$,
	\[
	\|Fx(t)-x_0\|\le \int_{[t_0,t)}\omega(\|x(s)-x_0\|)\dif \overline{g}(s)+\sum_{i=1}^n \int_{[t_0,\sigma)}|f_i(s,x_0)|\dif g_i(s)<R,
	\]
	for every $t\in \overline I_\sigma$. That is, $F(B)\subset B$. It remains to verify that $F(B)$ is relatively compact in $\mathcal{BC}_{\boldsymbol g}(\overline I_\sigma,\R^n)$. Firstly, note that for $x\in B$ and $i\in\{1,2,\dots,n\}$,
	\[
	|f_i(t,x(t))|\le |f_i(t,x(t))-f_i(t,x_0)|+|f_i(t,x_0)|\le \varphi(t)\omega(R)+|f_i(t,x_0)|=:M_i(t),
	\]
	for $g_i$--a.a. $t\in I_\sigma$. Observe that $F(B)\subset \mathcal{AC}_{\boldsymbol g}(\overline I_\sigma,\R^n)$ and
	\[(\left(Fx\right)_i)'_{g_i}(t)=f_i(t,x(t)),\quad g_i\mbox{--a.a. $t\in I_\sigma$,}\quad  \mbox{$x\in B$},\quad i=1,2,\dots,n.\]
	Since $M_i\in\mathcal L^1_{g_i}(I_\sigma,[0,+\infty))$ for each $i\in\{1,2,\dots, n\}$, it follows from Proposition~\ref{multirelcomp} that $F(B)$ is relatively compact in $\mathcal{BC}_{\boldsymbol g}(\overline I_\sigma,\R^n)$. Thus, Schauder's Fixed Point Theorem guarantees the existence of solution of \eqref{Multivp} on $I_\sigma$, while the uniqueness is a consequence of Theorem~\ref{multiosgood}.
\qed

As a final note, we provide an example where the results obtained in this paper can be applied while those in \cite{LoMa19} cannot. To that end, we consider the maps $\exp^{[k]}:\mathbb R\to\mathbb R$, $k=0,1,2,\dots,$ defined
\[\exp^{[0]}(t)=t,\quad \exp^{[k]}(t)=\exp(\exp^{[k-1]}(t)),\quad k\in\mathbb N,\]
the functions $\log^{[k]}:(\exp^{[k-1]}(0),+\infty)\to\mathbb R$, $k\in\mathbb N$, defined as
\[\log^{[1]}(t)=\log(t),\quad \log^{[k]}(t)=\log(\log^{[k-1]}(t)),\quad k\ge 2,\]
and the family of functions $\omega_k:[0,+\infty)\to\mathbb R$, $k\in\mathbb N$, defined as
\begin{equation}\label{eq-fn}
\omega_k(t)=\begin{cases}
0&\mbox{if }t=0,\\
\displaystyle t\prod_{j=1}^k \log^{[j]}\left(\frac{1}{t}\right)&\mbox{if } \displaystyle 0<t<\frac{1}{e_k},\\
\displaystyle \frac{1}{e_k^2}\prod_{j=1}^k \exp^{[j]}(1)&\mbox{if }\displaystyle t\ge\frac{1}{e_k},
\end{cases}
\end{equation}
where $e_k:=\exp^{[k]}(1)$. In \cite[Proposition~5.1]{MaMon20}, the authors proved some interesting properties of $\omega_k$, $k\in\mathbb N$. In particular, they showed that $\omega_k$ is nondecreasing, differentiable in $(0,1/e_k)$ and
\[\lim_{t\to 0^+}\omega_k'(t)=+\infty.\]
Furthermore, it is shown that $\omega_k$, $k\in\mathbb N$, satisfies \eqref{osgoodcond} for any $u_0>0$ and, in \cite[Example~5.2]{MaMon20}, that for any $x,y\in\mathbb R$,
\begin{equation}\label{subadditive}
|\omega_k(|x|)-\omega_k(|y|)|\le \omega_k(|x-y|).
\end{equation}
With this, we have all the information necessary for the following example.
\begin{example}
	Let $\boldsymbol g:\mathbb R\to\mathbb R^n$, $\boldsymbol g=(g_1,g_2,\dots, g_n)$, such that each $g_i$, $i\in\{1,2,\dots,n\}$, is a derivator; and $\varphi:[0,1)\to[0,+\infty)$ be a Borel measurable map which is bounded. Observe that this implies that $\varphi\in\mathcal L^1_{g_i}([0,1),[0,+\infty))$ for each $i\in\{1,2,\dots,n\}$. For each $k\in\mathbb N$, consider the initial value problem
	\begin{equation}\label{multiivp-n}
	x'_{\boldsymbol g}(t)=f_k(t,x(t)),\quad g\mbox{--a.a. }t\in[0,1),\quad x(0)=x_0,
	\end{equation}
	where $f_k:[0,1)\times\R^n \to \R^n$, $f_k=(f_{k,1},f_{k,2},\dots, f_{k,n})$, is given by
	\[f_{k,i}(t,x)=\varphi(t)\,\omega_k(\|x\|),\quad (t,x)\in\mathbb [0,1)\times\mathbb R,\quad i\in\{1,2,\dots,n\},\]
	with $\omega_k$ as in \eqref{eq-fn}. Note that each $f_{k,i}$, $i\in\{1,2,\dots,n\}$, $k\in\mathbb N$, does not satisfied a Lipschitz condition on the whole $\R$ as the derivative of $\omega_k$ is unbounded on any neighbourhood of $0$. Therefore, the hypotheses of \textup{\cite[Theorem~4.3]{LoMa19}} cannot be satisfied. Similarly, \textup{\cite[Theorem~4.4]{LoMa19}} cannot be applied for any ball around $x_0$ containing $0$.
	However,  \textup{Theorem~\ref{multiexist-Schauder}} yields that \eqref{multiivp-n} has a unique solution.
	Indeed, it is clear that conditions \textup{(i)} and \textup{(ii)} are satisfied. Furthermore, as mentioned before, \textup{\cite[Proposition~5.1]{MaMon20}} guarantees that condition \textup{(iii)} is satisfied. Now, for condition \textup{(iv)}, for each $i\in\{1,2,\dots,n\}$, the fact that $\omega_k$ is nondecreasing and \eqref{subadditive} yield
	\[|f_{k,i}(t,x)-f_{k,i}(t,y)|= \varphi(t)|\omega_k(\|x\|)-\omega_k(\|y\|)|\le \varphi(t)\omega_k(|\|x\|-\|y\||)\le \varphi(t)\omega_k(\|x-y\|),\quad k\in\mathbb N,\]
	which shows that condition \textup{(iv)} is satisfied. Hence, the hypotheses of \textup{Theorem~\ref{multiexist-Schauder}} are satisfied. Using analogous reasonings, we can show that the hypotheses of \textup{Theorem~\ref{multiexist-uniq}} are satisfied if we consider $f_k$ to be defined on $[0,1)\times \overline{B(x_0,r)}$, $r>0$, instead of on $[0,1)\times\mathbb R$.
\end{example}

\section*{Acknowledgments}
Ignacio M\'arquez Alb\'es was partially supported by Xunta de Galicia under grant ED481A-2017/095 and project ED431C 2019/02.

F. Adri\'an F. Tojo was partially supported by Xunta de Galicia, project ED431C 2019/02, and by the Agencia Estatal de Investigaci\'on (AEI) of Spain under grant MTM2016-75140-P, co-financed by the European Community fund FEDER.

\end{document}